\pdfoutput=1  

\documentclass[11pt]{article}

\usepackage[margin=0.8in]{geometry}
\usepackage{subcaption}
\usepackage{adjustbox}
\usepackage{tikz}
\usetikzlibrary{shapes.callouts,arrows.meta}
\usepackage{algorithm}
\usepackage{algpseudocode}

\usepackage[most]{tcolorbox}
\usepackage{xcolor}
\usepackage{booktabs,array,multirow,longtable}
\newtcolorbox{callout}[2][]{
    enhanced,
    breakable,
    colback=blue!5,
    colframe=black!10,
    coltitle=black,
    fonttitle=\bfseries,
    title=#2,
    arc=2mm,
    boxrule=0.8pt,
    left=3mm,
    right=3mm,
    top=2mm,
    bottom=2mm,
    before skip=8pt,
    after skip=8pt,
    #1
}
\usepackage{amsmath,xparse,etoolbox}
\usetikzlibrary{calc,decorations.markings,shadows.blur}

\def\VertexRadius{0.23cm}
\def\FingerWidth{2pt}
\def\PropWidth{3.5pt}
\def\CrossWidth{1.2pt}
\def\CrossSize{0.1cm}

\colorlet{VertexColor}{gray!72}
\colorlet{FingerColor}{blue!65}
\colorlet{PropColor}{red!78}
\colorlet{CrossColor}{black!98}

\newcommand{\Tube}[3]{%
  \begingroup
  \pgfmathsetlengthmacro{\ShadowW}{#2+2pt}
  \pgfmathsetlengthmacro{\RimW}{#2+1.4pt}
  \pgfmathsetlengthmacro{\ShadeW}{#2/3}
  \pgfmathsetlengthmacro{\LightW}{#2/7}
  \pgfmathsetlengthmacro{\TubeShift}{#2/4}
 \draw[black,opacity=.2,line width=\ShadowW, line cap=round,line join=round, transform canvas={xshift=.7pt,yshift=-1.2pt}] #3;
  \draw[#1!72!black,line width=\RimW,       line cap=round,line join=round] #3;
  \draw[#1,line width=#2,line cap=round,line join=round] #3;
  \endgroup
}

\newcommand{\CrossMark}{%
  \draw[CrossColor!70!black,line width=2pt,line cap=round]
    (-\CrossSize,-\CrossSize)--(\CrossSize,\CrossSize)
    (-\CrossSize, \CrossSize)--(\CrossSize,-\CrossSize);
  \draw[CrossColor,line width=\CrossWidth,line cap=round]
    (-\CrossSize,-\CrossSize)--(\CrossSize,\CrossSize)
    (-\CrossSize, \CrossSize)--(\CrossSize,-\CrossSize);
  \draw[white,opacity=.45,line width=.45pt,line cap=round]
    (-.11cm,.11cm)--(.11cm,-.11cm);
}

\newcommand{\PutCrossings}[2]{%
  \def\temp{#2}%
  \ifx\temp\empty\else
    \foreach \t in {#2}{%
      \path[postaction={decorate},decoration={markings,
        mark=at position \t with {\CrossMark}}] #1;
    }%
  \fi
}

\newcommand{\VertexList}{}

\NewDocumentCommand{\D}{O{1}+m}{%
  \begin{tikzpicture}[x=1cm,y=1cm,scale=#1,transform shape]
    \def\VertexList{}%
    #2
    \VertexList
  \end{tikzpicture}%
}

\NewDocumentCommand{\V}{m m O{}}{%
  \coordinate (#1) at (#2);
  \appto\VertexList{%
    \shade[ball color=VertexColor,draw=black!68!black,line width=.5pt,
           blur shadow={shadow opacity=.20,shadow blur steps=5,shadow xshift=.8pt,shadow yshift=-1.2pt}]
      (#1) circle[radius=\VertexRadius];
    \fill[white,opacity=.92] ($(#1)+(-.12,.13)$) circle[radius=.045cm];
    \if\relax\detokenize{#3}\relax\else
      \node[font=\large] at (#1) {$#3$};
    \fi
  }%
}

\NewDocumentCommand{\E}{m m O{} O{}}{%
  \def\Path{(#1)--(#2)}%
  \Tube{PropColor}{\PropWidth}{\Path}
  \PutCrossings{\Path}{#4}
  \if\relax\detokenize{#3}\relax\else
    \path[postaction={decorate},decoration={markings,
      mark=at position .5 with {\node[font=\large,above=7pt] {$#3$};}}] \Path;
  \fi
}

\newcommand{\F}[6]{%
  \def\Path{(#1)--($(#1)+(#2:#3)$)}%

  \Tube{FingerColor}{\FingerWidth}{\Path}

  \PutCrossings{\Path}{#5}

\foreach \t [count=\i] in {#5}{%
  \foreach \lab [count=\j] in {#6}{%
    \ifnum\i=\j
      \path[
        postaction={decorate},
        decoration={
          markings,
          mark=at position \t with {
            \pgftransformresetnontranslations
            \pgftransformyshift{.1mm}
            \node[
              font=\scriptsize,
              fill=none,
              fill opacity=.9,
              text opacity=1,
              inner sep=.6pt
            ] at (4.5mm,0) {$\lab$};
          }
        }
      ] \Path;
    \fi
  }%
}%

\node[
  font=\large,
  inner sep=1pt
] at ($(#1)+(#2:#3+.35)$) {$#4$};
}
\NewDocumentCommand{\FC}{m m m m O{}}{%
  \pgfmathsetmacro{\outangle}{#2+#4}
  \pgfmathsetmacro{\inangle}{#2+180-#4}
  \def\Path{(#1) to[out=\outangle,in=\inangle,looseness=1.04] ($(#1)+(#2:#3)$)}%
  \Tube{FingerColor}{\FingerWidth}{\Path}
  \PutCrossings{\Path}{#5}
}

\usepackage{graphicx}
\DeclareGraphicsExtensions{.pdf,.png,.jpg}

\usepackage{amsmath,amssymb,amsthm,mathtools,bm}
\usepackage{braket}
\usepackage[colorlinks=true,citecolor=blue]{hyperref}
\usepackage{cite}

\usepackage{stmaryrd}    
\usepackage{xfrac}
\DeclareMathAlphabet{\mathpzc}{OT1}{pzc}{m}{it} 
\newcommand\mathscr[1]{\scalebox{1.1}{$\mathpzc{#1}$}}
\usepackage{booktabs} 

\usepackage[all]{xy}
\usepackage{tikz}
\usetikzlibrary{cd}
\usetikzlibrary{calc}

\usepackage{tikz}
\usetikzlibrary{arrows.meta}
\usepackage{float}
\usepackage{booktabs} 
\usetikzlibrary{decorations.markings, arrows.meta, calc, shapes.misc}
\usepackage{setspace}
\usepackage{float}
\usepackage{xcolor}
\usepackage{rotating}
\usepackage{multirow}
\usepackage{colortbl}
\usepackage{amsthm}
\usepackage{amssymb}

\usepackage{hyperref}
\hypersetup{
  colorlinks=true,
  allcolors=black
}
\usepackage{cleveref}
\crefformat{section}{\S#2#1#3}
\crefformat{subsection}{\S#2#1#3}
\crefformat{subsubsection}{\S#2#1#3}

\usepackage{titlesec}
\titleformat{\section}{\normalfont\sffamily\large\bfseries}{\thesection.}{0.5em}{}
\titleformat{\subsection}{\normalfont\sffamily\bfseries}{}{0em}{}
\titlespacing*{\section}{0pt}{1.4ex plus .2ex}{0.6ex}
\titlespacing*{\subsection}{0pt}{1.0ex plus .2ex}{0.4ex}

\usepackage{mathptmx}
\usepackage[cal=cm]{mathalfa}  

\makeatletter
\newif\if@sup
\newtoks\@sups
\def\append@sup#1{\edef\act{\noexpand\@sups={\the\@sups #1}}\act}%
\def\reset@sup{\@supfalse\@sups={}}%
\def\mk@scripts#1#2{\if #2/ \if@sup ^{\the\@sups}\fi \else%
  \ifx #1_ \if@sup ^{\the\@sups}\reset@sup \fi {}_{#2}%
  \else \append@sup#2 \@suptrue \fi%
  \expandafter\mk@scripts\fi}
\def\tensor#1#2{\reset@sup#1\mk@scripts#2_/}
\def\multiscripts#1#2#3{\reset@sup{}\mk@scripts#1_/#2%
  \reset@sup\mk@scripts#3_/}
\makeatother

\makeatletter
\newbox\slashbox \setbox\slashbox=\hbox{$/$}
\def\itex@pslash#1{\setbox\@tempboxa=\hbox{$#1$}
  \@tempdima=0.5\wd\slashbox \advance\@tempdima 0.5\wd\@tempboxa
  \copy\slashbox \kern-\@tempdima \box\@tempboxa}
\def\slash{\protect\itex@pslash}
\makeatother

\usepackage[titletoc]{appendix}
\theoremstyle{italics}
\newtheorem{theorem}{Theorem}[section]

\newtheorem{prop}[theorem]{Proposition}

\theoremstyle{defn}

\let\PLAINthebibliography\thebibliography
\renewcommand\thebibliography[1]{
  \PLAINthebibliography{#1}
  \setlength{\parskip}{0.5pt}
  \setlength{\itemsep}{0.5pt plus .3ex}
}

\title{A TQFT-based Platform for Efficient Computation of Knot Invariants  
}

\author{
\def\arraystretch{.9}
  \begin{tabular}{c}
  Amena Al Rawi\rlap{${}^{a, \, c}$}
  \\
  {\color{gray}\footnotesize\tt aa12439@nyu.edu}
  \end{tabular} 
  \quad 
  \def\arraystretch{.9}
  \begin{tabular}{c}
  Hisham Sati\rlap{${}^{{b},\,{c}, \, d}$}
  \\
  {\color{gray}\footnotesize \tt hsati@nyu.edu}
  \end{tabular}
  \;\;\;\;
  \def\arraystretch{.9}
  \begin{tabular}{c}
  Vivek Kumar Singh\rlap{${}^{\,c}$}
  \\
  {\color{gray}\footnotesize\tt vks2024@nyu.edu}
  \end{tabular}
}

\date{}
\begin{document}
\maketitle
\vspace{2pt}
\vspace{8pt}
\begin{abstract}
 We present an interactive web platform\cite{RSS} that unifies the construction Feynman ribbon diagrams (FRDs), the evaluation of higher-rank Chern--Simons knot invariants, and the identification of FRD-like knots at higher crossing numbers. These tree-structured diagrams naturally represent arborescent knots, which we refer to throughout as \emph{FRD-like knots}. Within a single visual environment, users can construct an FRD as a tensor network, evaluate its associated Chern--Simons invariants, and use the resulting invariant data to distinguish and identify the corresponding knot. To our knowledge, this is the first platform to combine diagrammatic construction, tensor-network evaluation, invariant computation, and knot identification within a unified workflow.
\end{abstract} 

\vspace{1cm} 
\begin{center}
\begin{minipage}{12cm}
\tableofcontents
\end{minipage}
\end{center}




\vfill

\hrule
\vspace{-1pt}

{
\hypertarget{DoS}{}
\footnotesize
\noindent
\def\arraystretch{1}
\tabcolsep=0pt

\begin{tabular}{ll}
${}^a$\,
&
Computer Science, Division of Science,
New York University Abu Dhabi, UAE.  
\\
${}^b$\,
&
Mathematics, Division of Science, 
New York University Abu Dhabi, UAE.  
\\
${}^c$ & 
Center for Quantum and Topological Systems (CQTS),
\\
&
NYUAD Research Institute,
New York University Abu Dhabi, UAE.  
\\
${}^d$\, & The Courant Institute for Mathematical Sciences, NYU, NY.
\end{tabular}
\hfill
\adjustbox{raise=-5pt}{
\href{https://ncatlab.org/nlab/show/Center+for+Quantum+and+Topological+Systems}{\includegraphics[width=3cm]{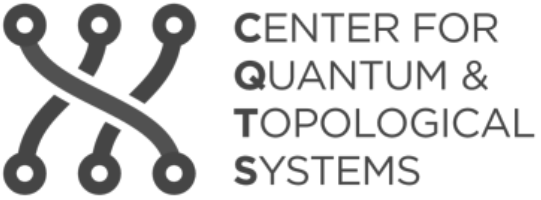}}
}

\vspace{0.1mm}

\noindent
The authors acknowledge the support by {\it Tamkeen} under the 
{\it NYU Abu Dhabi Research Institute grant} {\tt CG008}.
}

\section{Introduction}

The central problem of knot theory is to determine when two closed curves in three-dimensional space represent the same knot. Knot invariants make this problem more tractable, but the most informative invariants, especially colored HOMFLY--PT polynomials \cite{HOMFLY1985, PT1987}, can be difficult to compute for large crossing knots and higher representations.

The present work focuses on arborescent (double-fat) knots that admit a description in terms of \emph{Feynman ribbon diagrams} (FRDs). The structured, tree-like form of an FRD can be translated into a tensor network whose components are evaluated using Topological quantum Field theory (TQFT). We developed the \textbf{TQFT Knot Explorer}\cite{RSS} to facilitate this workflow. The tool enables users to draw an FRD, construct and evaluate the associated tensor network, compute the corresponding knot invariants, and compare the results with entries in a classification database.

Benchmarks against the \texttt{KnotTheory} colored Jones algorithm and the braid-walk method~\cite{KnotAtlas, HajijLevitt2018} show that our approach is faster in most tested cases, establishing it as an efficient and scalable tool for knot classification and invariant computation.

Converting technically demanding calculations into intuitive, reproducible diagrammatics allows the platform to substantially advance the study of higher-rank Chern--Simons invariants and broaden their accessibility. A central result is a remarkably compact two-vertex FRD whose invariant data completely distinguishes knots with up to $10$ crossings and continues to enable the effective identification of FRD-like knots at higher crossing numbers. This demonstrates that even a minimal FRD architecture can encode unexpectedly rich topological information.

\paragraph{Knots and the Classification Problem.}

A knot is an embedding of the circle \(S^1\) in \(\mathbb{R}^3\), or equivalently in \(S^3\), considered up to continuous deformation without cutting or gluing. A link is a finite collection of disjoint knots. Knots are usually represented by planar diagrams containing overcrossings and undercrossings. By Reidemeister's theorem, two diagrams represent the same knot if and only if they are related by planar deformations and a finite sequence of Reidemeister moves~\cite{Reidemeister1927}. Although this gives a theoretical test for equivalence, finding such a sequence can be extremely difficult in practice. Classification becomes even harder as the crossing number increases. The number of prime knots grows from \(165\) at \(10\) crossings to more than \(1.7\) million at \(16\) crossings~\cite{HosteThistlethwaiteWeeks}. Moreover, classification is not merely the generation of diagrams: many different diagrams can represent the same knot, and one must also determine whether a diagram has the smallest possible number of crossings.

\paragraph{From Classical to Polynomial Invariants.}
To avoid direct diagram-by-diagram comparison, knot theorists use \emph{knot invariants}: quantities that remain unchanged under Reidemeister moves. If two knots have different invariant values, then they are certainly inequivalent. Equal values, however, do not always imply that the knots are the same. The Alexander polynomial was one of the first effective polynomial invariants~\cite{Alexander1928}. It distinguishes many knots, though it generally fails to distinguish a knot from its mirror image. The Jones polynomial introduced a stronger invariant and revealed important connections between knot theory and mathematical physics~\cite{Jones1985}. The two-variable HOMFLY--PT polynomial \cite{HOMFLY1985, PT1987} further unifies and extends these constructions, reducing to the Jones and Alexander polynomials under appropriate specializations. These polynomials greatly improve knot classification, but none of them is a complete invariant. Distinct knots may share the same Alexander, Jones, or HOMFLY--PT polynomial. Classification therefore remains difficult even when polynomial information is available.

\paragraph{Colored Invariants and the Computational Bottleneck.}
A deeper interpretation of polynomial knot invariants emerges within the framework of topological quantum field theory (TQFT). In particular, Chern--Simons theory realizes these invariants as expectation values of Wilson-loop operators~\cite{WitCS}\footnote{The corresponding algebraic construction is the
Reshetikhin--Turaev formalism, which uses quantum-group representations\cite{RT,RT3M}.}.
Direct evaluation of Chern--Simons invariants is generally difficult because the path integral involves an infinite-dimensional integration over gauge fields (cf. \cite{FG91}). To make the computation tractable, we use the correspondence between Chern--Simons theory, with group $G$ and level $k$, on a three-manifold with a punctured $\mathbf{S}^2$ boundary and the conformal blocks of the $\widehat{G}_k$ Wess--Zumino--Novikov--Witten (WZNW) model, where $\widehat{G}_k$ is the associated Kac-Moody group of $G$. Under this correspondence, the Chern--Simons path integral on a three-ball defines a state in the boundary Hilbert space, identified with the corresponding space of conformal blocks \cite{WitCS}. The braiding of $\widehat{G}_k$ WZNW conformal blocks provides a systematic construction of knot and link invariants in which each component is labeled by an arbitrary representation $R$ of $G$~\cite{RGK1993}.

\medskip 
For $G=SU(N)$, assigning a representation $R$ to a knot defines its \emph{$R$-colored Chern--Simons invariant}, equivalently encoded by the corresponding \emph{$R$-colored HOMFLY--PT polynomial}. In the special case $N=2$, this construction reduces to the \emph{colored Jones polynomial}. Colored polynomials contain substantially
more information than the fundamental polynomial and are important in the
study of mutation, recursion relations, differential expansions, and
special-polynomial limits~\cite{DGR,evo,MMSle,Mort,Morton}. They also
appear naturally in topological string
theory~\cite{GopakumarVafa,OoguriVafa}.

Their computational cost, however, grows rapidly with the representation and
the complexity of the knot. Higher representations require larger Racah
matrices, more braiding data, and increasingly complicated tensor
contractions. This creates a practical gap between the theoretical strength of
colored invariants and their accessibility to researchers.

\section{TQFT Knot Explorer}
\label{sec:knot-diagram-Explorer}

Knots are most naturally understood through diagrams, whereas their quantum
invariants are often presented through algebraic formulas that tend to obscure the
underlying geometric structure. The \textbf{TQFT Knot Explorer}\cite{RSS} brings
these two perspectives together in a browser-based platform. Users can assemble
a knot from a small collection of graphical building blocks, assign integer
twist parameters and a representation label, and compute the corresponding
quantum invariant. Thus, every step from the knot diagram to the exact
calculation remains visible, reproducible, and readily comparable across
different examples.

\begin{figure}[ht]
  \centering
  \hspace{-1.1cm}
  \begin{subfigure}[c]{0.25\textwidth}
    \centering
    \includegraphics[width=\linewidth]{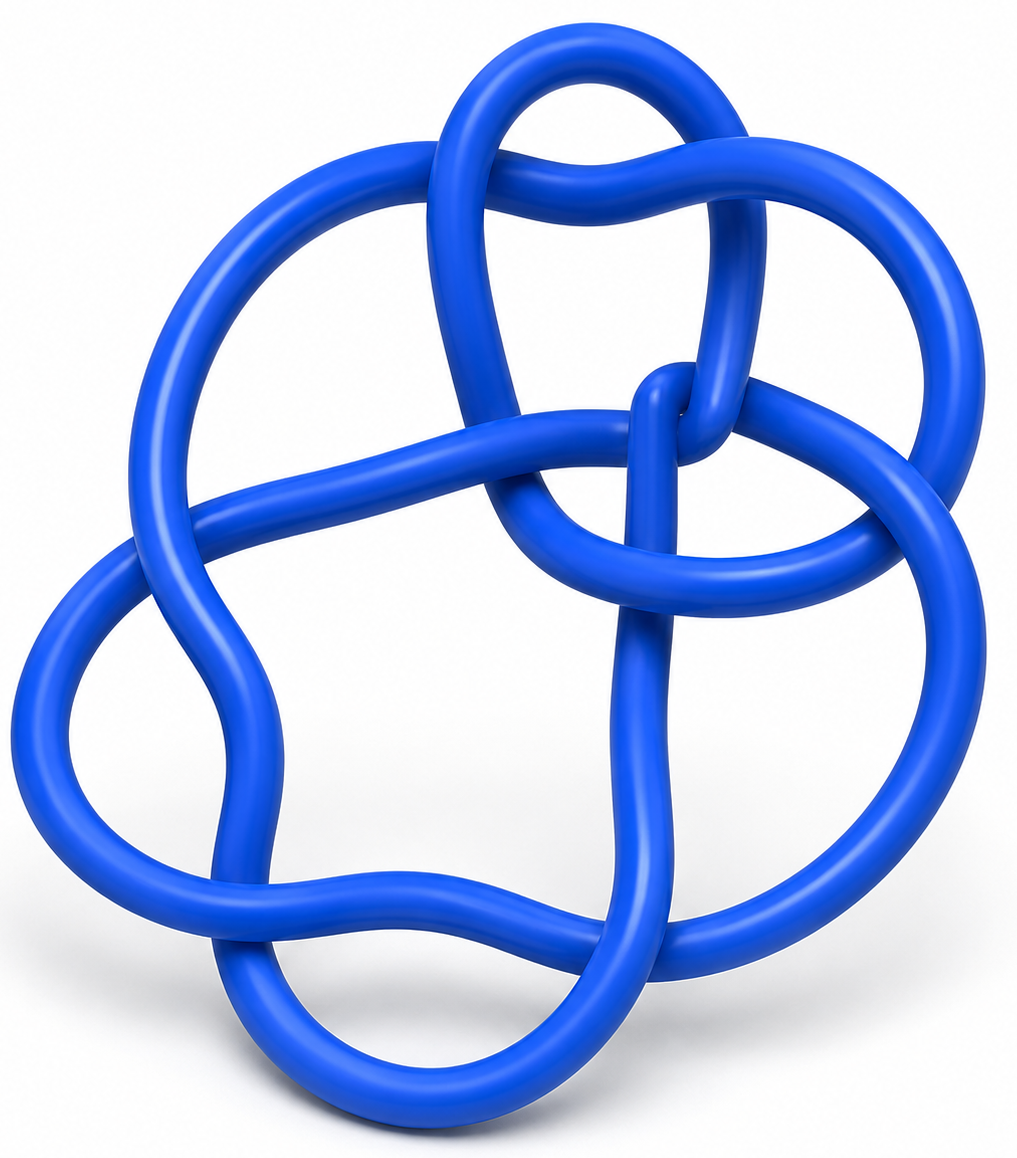}
  \end{subfigure}
  \hspace{0.005\textwidth}
  $=$
  \hspace{0.005\textwidth}
  \begin{subfigure}[c]{0.36\textwidth}
    \centering
    \includegraphics[width=\linewidth]{FD193.pdf}
  \end{subfigure}
  \hspace{0.005\textwidth}
  $\equiv$
  \hspace{0.005\textwidth}
  \begin{subfigure}[c]{0.25\textwidth}
    \centering
    \D[1]{
  \V{A}{-4,0}
  \V{B}{ -1,0}
  \E{A}{B}[1][.5]
  \F{A}{120}{2.2}{}{.5}{2~~}
  \F{A}{-120}{2.2}{}{.5}{3~~}
  \F{B}{65}{2.2}{}{.5}{3}
  \F{B}{-65}{2.2}{}{.5}{4}
}
\end{subfigure}
\vspace{-3mm} 
  \caption{Three complementary views of the knot $10_{93}$: its conventional
  projection, its FRD description, and the corresponding TQFT tangle
  representation.}
  \label{fig:knot-10-93}
\end{figure}

The current version of the Explorer focuses on \emph{arborescent knots}, often
described in this context as \emph{double-fat knots}. This extensive family
contains several familiar classes, including pretzel knots and $2$-bridge
knots; further discussions of these families may be found in
~\cite{Con,Caudron,BS}. A complete classification in terms of
such diagrammatic presentations becomes increasingly intricate at higher
crossing numbers. As an illustrative example, we show in Fig. \ref{fig:knot-10-93} how the knot $10_{93}$
can be represented as a double-fat knot. 

Its equivalent description arises
from the braiding and fusion operations acting on four-point conformal blocks
of the $\widehat{G}_{k}$ WZNW model.

The Explorer encodes these constructions using \emph{Feynman ribbon diagrams},
abbreviated as FRDs, and we use the term \emph{FRD-like knots} for knots
presented in this form. This terminology refers specifically to the
representation implemented in the Explorer and does not imply that every knot
admits a tree-shaped FRD description.

\begin{figure}[ht]
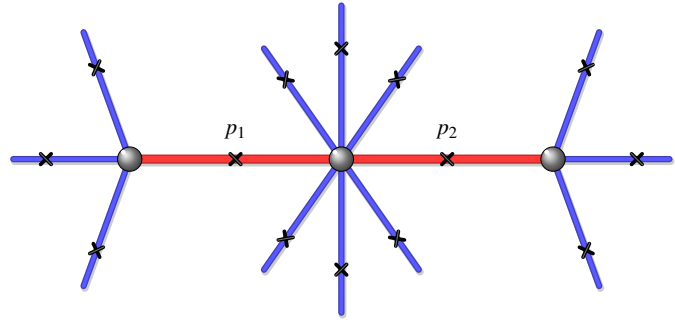

\centering
\adjustbox{scale=0.7}{
\D[1]{
  \V{A}{-4,0}
  \V{B}{ 0,0}
  \V{C}{ 4,0}

  \E{A}{B}[p_1][.5]
  \E{B}{C}[p_2][.5]

  \F{A}{180}{2.2}{}{.72}{}
  \F{A}{110}{2.55}{}{.72}{}
  \F{A}{-110}{2.55}{}{.72}{}

  \F{B}{90}{2.9}{}{.72}{}
  \F{B}{55}{2.55}{}{.72}{}
  \F{B}{125}{2.55}{}{.72}{}
  \F{B}{-55}{2.55}{}{.72}{}
  \F{B}{-90}{2.9}{}{.72}{}
  \F{B}{-125}{2.55}{}{.72}{}

  \F{C}{0}{2.2}{}{.72}{}
  \F{C}{70}{2.55}{}{.72}{}
  \F{C}{-70}{2.55}{}{.72}{}
}
}
\vspace{-4mm} 
\caption{The basic structure of a FRD. Blue lines are
fingers, red lines are propagators, and gray dots are vertices. Small black marks
on the fingers indicate positions carrying integer-valued braiding
parameters.}
\label{fig:FRD}
\end{figure}

\paragraph{A visual language for \emph{FRD-like knots}.}
An FRD may be interpreted as a graphical construction plan for the corresponding knot, as illustrated in Fig.~\ref{fig:FRD}.
\vspace{-2mm} 
The gray dots are \emph{vertices},
where different parts of the diagram meet. A blue line attached to a single vertex
is a \emph{finger}, also called a \emph{leg}, while a line joining two vertices
is a \emph{propagator}. Small marks placed along a finger indicate braiding
positions. Each mark carries a signed integer that controls the amount and
direction of local twisting. A propagator likewise carries an integer that
records the twist between the two vertices it connects.

\smallskip 
These few ingredients separate the architecture of the knot from its local
twisting. The vertices, fingers, and propagators determine how the diagram is
assembled; the integers determine how strongly its different regions are
twisted. The same underlying FRD can therefore describe an entire family of
knots as its braiding parameters vary.
 
\medskip 
To enable a systematic construction of knots within the FRD framework, we
introduce the fundamental graphical building blocks and their TQFT analogues,
leading to the summary in Table~\ref{tqft-tensor} below. General FRDs can then be
assembled from elementary vertices, fingers, and propagators.

\begin{figure}[ht]
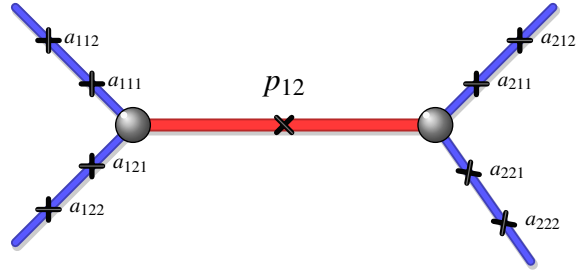

\centering
\D[1]{
  \V{A}{-4,0}
  \V{B}{ 0,0}
  \E{A}{B}[p_{12}][.5]
  \F{A}{135}{2.2}{}{.35,.72}{a_{111},a_{112}}
  \F{A}{-135}{2.2}{}{.35,.72}{~~a_{121},~~a_{122}}
  \F{B}{45}{2.2}{}{.35,.72}{~~a_{211},~~a_{212}}
  \F{B}{-55}{2.2}{}{.35,.72}{~~a_{221},~~a_{222}}
}
\vspace{-4mm} 
\caption{A two-vertex Feynman ribbon diagram with four fingers and a connecting
propagator.}
\label{fig:FRD2}
\end{figure}

We introduce the two-vertex Feynman ribbon diagram shown in
Fig.~\ref{fig:FRD2}, denoted by $\Gamma(A_1;p_{12};A_2)$. It consists of two vertices connected by a propagator \(p_{12}\), with two
fingers attached to each vertex.

The data associated with the two vertices are
\[
A_1=
\big[(a_{111},a_{112}),(a_{121},a_{122})\big],
\qquad
A_2=
\big[(a_{211},a_{212}),(a_{221},a_{222})\big],
\]
where each ordered pair specifies the two twist parameters associated with a
single finger. The label \(p_{12}\) determines the propagator connecting the
two vertices. Such Feynman ribbon diagrams are also referred to as evolution
families (\cite{MMM13,MMMRSS2017}). Our main finding here is the following (which we will expand on in \cite{RSS2026}):

\begin{prop}
The evolution family of two-vertex Feynman ribbon diagrams
\(\Gamma(A_1;p_{12};A_2)\) is complete, within the class of FRD-like knots,
through ten crossings. More precisely, every FRD-like knot \(K\) with crossing
number \(c(K)\leq 10\) admits a representation of the form
\(\Gamma(A_1;p_{12};A_2)\) for an appropriate choice of the finger parameters
\(A_1,A_2\) and the propagator \(p_{12}\).
\end{prop}
The complete classification of knots up to $10$ crossings is presented in Appendix \ref{Appendix1}. Moreover, the diagram \(\Gamma(A_1; p_{12}; A_2)\) remains effective for identifying and representing many knots with higher crossing numbers, including numerous examples with up to $13$ crossings reported in \cite{RSS2026}.

\subsection{From a diagram to a quantum invariant}
\label{sec:tqft-building-blocks}

The real power of an FRD does not lie only in its simplified graphical representation. Each of its fundamental components has a corresponding building block in topological quantum field theory (TQFT). A vertex combines the information carried by its incident branches, a finger encodes the braiding data associated with a local twist region, and a propagator transfers information between two vertices. 

\smallskip 
Connecting these TQFT building blocks according to the combinatorial structure of the graphical diagram allows the Explorer to construct the tensor network associated with the knot. Using the amplitudes of these fundamental building blocks, described in detail in Ref.~\cite{RSS2026}, one can compute the full amplitude associated with any given tree-shaped FRD. The amplitudes of the individual blocks are derived from the conformal blocks, braiding operators, and fusion data of the corresponding four-point ${\rm SU}(N)$ WZNW conformal field theory. The complete derivation, together with the conventions and explicit tensor formulas, will be presented in a companion work\cite{RSS2026}. For the website reader, the essential correspondence is summarized in Table~\ref{tqft-tensor}.

\smallskip 
Here, \(\mathcal{H}\) denotes the space of conformal blocks of the underlying WZNW model. A user of the Explorer does not need to manipulate this space directly. Rather, it is the internal mathematical space in which the local information carried by vertices, fingers, and propagators is combined. Once the building blocks have been assembled, the Explorer connects their tensors in precisely the same pattern in which the corresponding graphical components are connected in the FRD. In mathematical terms, the internal indices are contracted. Schematically, the resulting amplitude may be written as
\[
\mathcal{A}_{R}(\Gamma)
=\operatorname{Contract}
\big(\{\mathcal{O}_{v}\},\{F_{f}\},\{P_{e}\}
\big)\,,
\]
where \(v\), \(f\), and \(e\) label the vertices, fingers, and propagator edges of the diagram, respectively. The tensors \(\mathcal{O}_{v}\), \(F_{f}\), and \(P_{e}\) represent the amplitudes assigned
to these fundamental components, respectively. In Chern--Simons theory, this amplitude is related to the expectation value of a Wilson loop \(W_R[K]\), colored by a representation \(R\):
\[
\mathcal{A}_{R}(\Gamma)
\propto
\bigl\langle W_R[K]\bigr\rangle \,.
\]
The precise proportionality factor depends on the normalization and framing conventions. After the appropriate normalization, the result is returned, in the platform's conventions, as a Laurent polynomial in the variables \(a\) and \(q\). This polynomial is the quantum invariant associated with the selected diagram and representation color.

Behind the interface, the symbolic engine evaluates the required braiding eigenvalues and Racah matrices \cite{Nawata2013,Alekseev2020,GJ,mmmrs,MM2016}.
The visual diagram is therefore not separate from the mathematics. It tells the engine which TQFT building blocks must be used, which parameters must be assigned to them, and how their tensor indices must be contracted. The entire procedure may be summarized schematically as
\begin{center} 
\colorbox{lightgray}{
$\text{Construct an FRD}
\;\longrightarrow\;
\text{evaluate the amplitude}
\;\longrightarrow\;
\text{knot identification} \;\longrightarrow\;
\text{colored knot polynomials}
$}
\end{center}
At present, the Explorer provides an efficient framework for computing colored quantum invariants of knots in SU$(2)$- Chern--Simons theory. Its computational engine accepts arbitrary admissible FRD representations and is not intrinsically restricted to knots with at most $13$ crossings. The present $13$-crossing limitation arises solely from the knot-identification database, which currently recognizes knots only within this range.  In future updates\cite{RSS}, the platform will support computations for knots with higher crossing numbers and a wider class of gauge groups, including SU$(N)$.

\begin{table}[ht]
\centering
\renewcommand{\arraystretch}{1.1}

\begin{tabular}{
|>{\centering\arraybackslash}m{0.13\textwidth}
|>{\centering\arraybackslash}m{0.18\textwidth}
|>{\centering\arraybackslash}m{0.23\textwidth}
|>{\centering\arraybackslash}m{0.36\textwidth}|
}
\hline
\rowcolor{lightgray}
\textbf{Element}
& \textbf{FRD block}
& \textbf{TQFT picture}
& \textbf{TQFT Amplitudes} \\
\hline

Vertex
&\adjustbox{valign=c}{
  \D[1]{
  \V{A}{-6,0}
  \F{A}{0}{1.2}{}{}{}
    \F{A}{45}{1.2}{}{}{}
    \F{A}{90}{1.2}{}{}{}
        \F{A}{135}{1.2}{}{}{}
   \F{A}{180}{1.2}{}{}{}
      \F{A}{-90}{1.2}{}{}{}
            \F{A}{-135}{1.2}{}{}{}

    \F{A}{-45}{1.2}{}{}{}

}
}
&
\adjustbox{valign=c}{
  \includegraphics[width=0.20\textwidth]{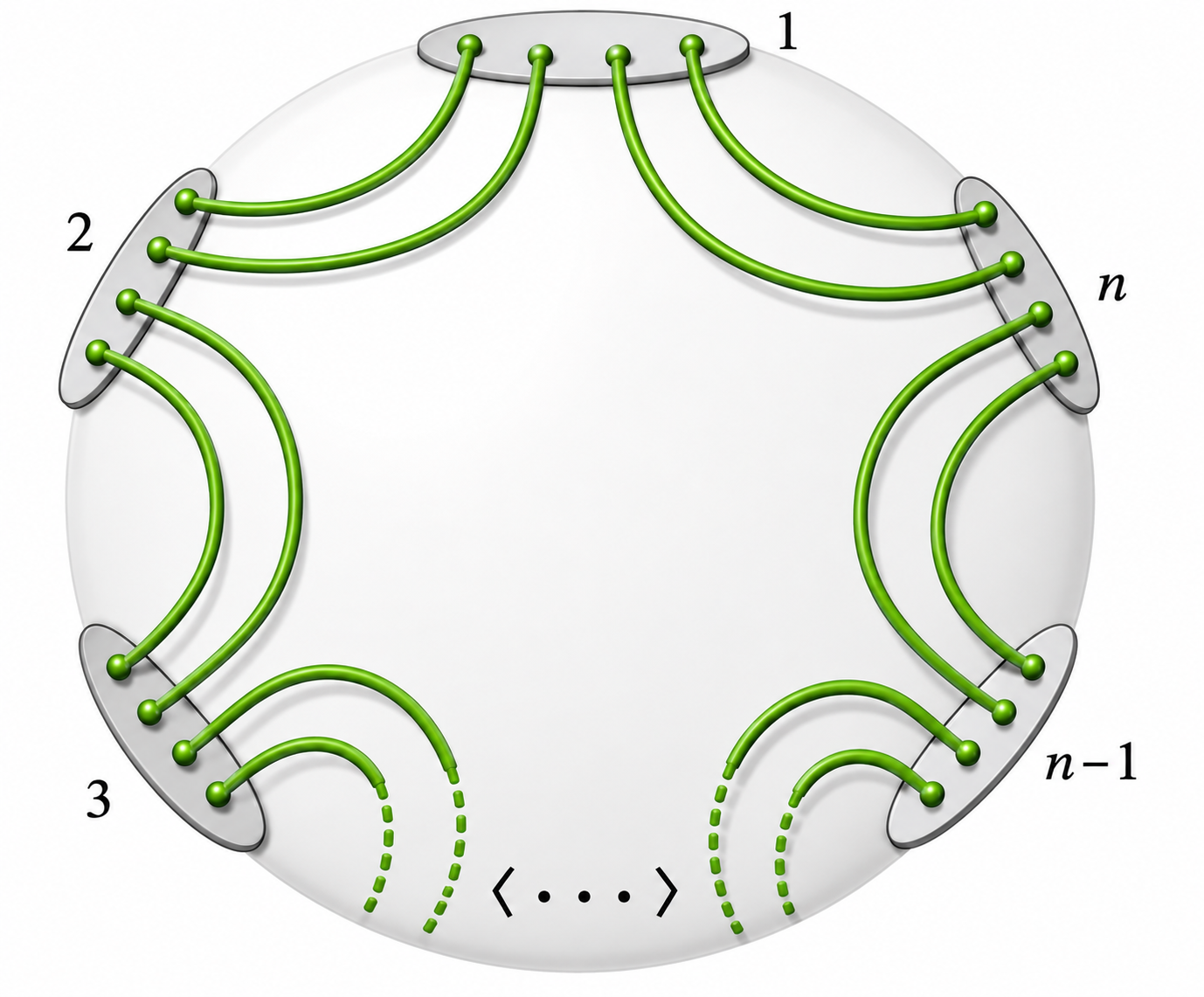}
}
&
The $n$-valent vertex($v$) amplitude:
\[
\mathcal{O}_{v}\in\mathcal{H}^{\otimes n}.
\]
\\
\hline

\multirow{2}{=}{\centering Fingers/Legs}
&
\adjustbox{valign=c}{
     \D[1]{
  \V{A}{-4,0}
  \F{A}{180}{2.2}{}{.5}{\hspace{-20pt}\raisebox{15pt}{$m$}}
}
}
&
\adjustbox{valign=c}{
  \includegraphics[width=0.20\textwidth]{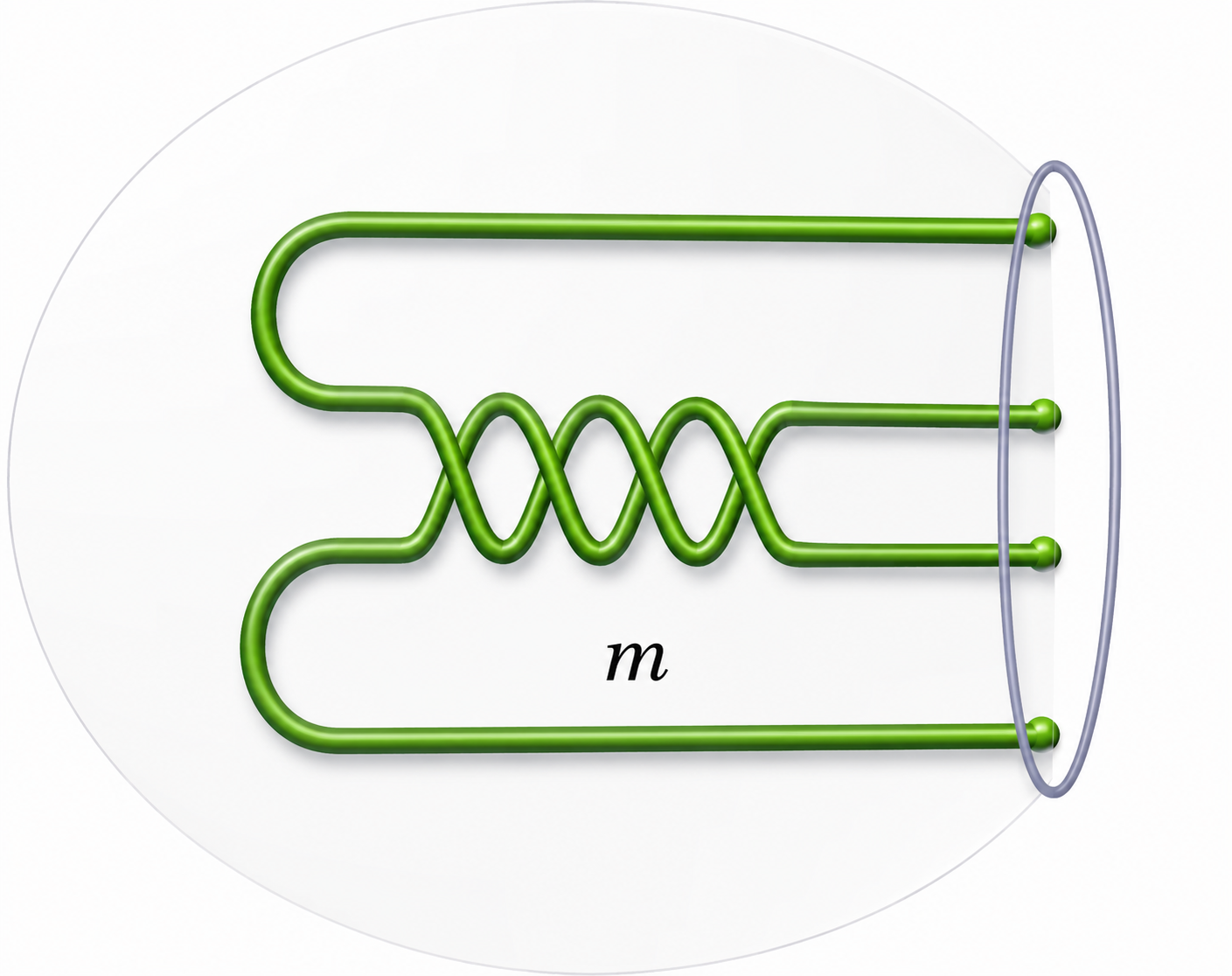}
}
&
The amplitude contributed by the legs connected to the $v$:
\[
F^{I}_{v}\in\mathcal{H}.
\]
\\
\cline{2-4}

&
\adjustbox{valign=c}{
       \D[1]{
  \V{A}{-4,0}
  \F{A}{180}{2.2}{}{.3,.7}{\hspace{-20pt}\raisebox{15pt}{$m$},\hspace{-20pt}\raisebox{15pt}{$n$}}
}
}
&
\adjustbox{valign=c}{
  \includegraphics[width=0.20\textwidth]{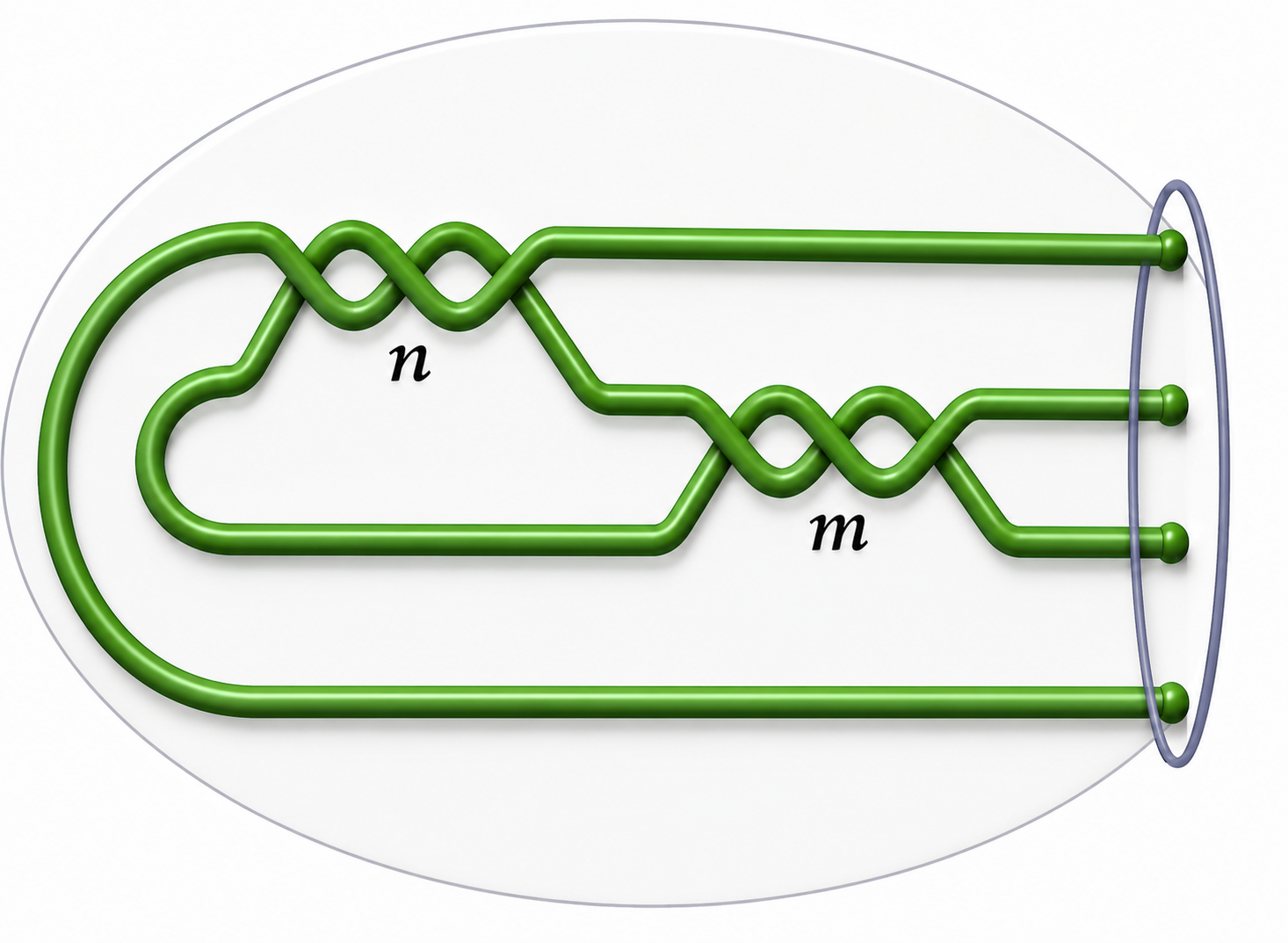}
}
&
The amplitude arising from the legs attached to vertex $v$:
\[
F^{II}_{v}\in\mathcal{H}.
\] 
\\
\hline

Propagator
&
\adjustbox{valign=c}{
       \D[1]{
  \V{A}{-2.5,0}
    \V{B}{0,0}
  \E{A}{B}[n][.5]
}
}
&
\adjustbox{valign=c}{
  \includegraphics[width=0.20\textwidth]{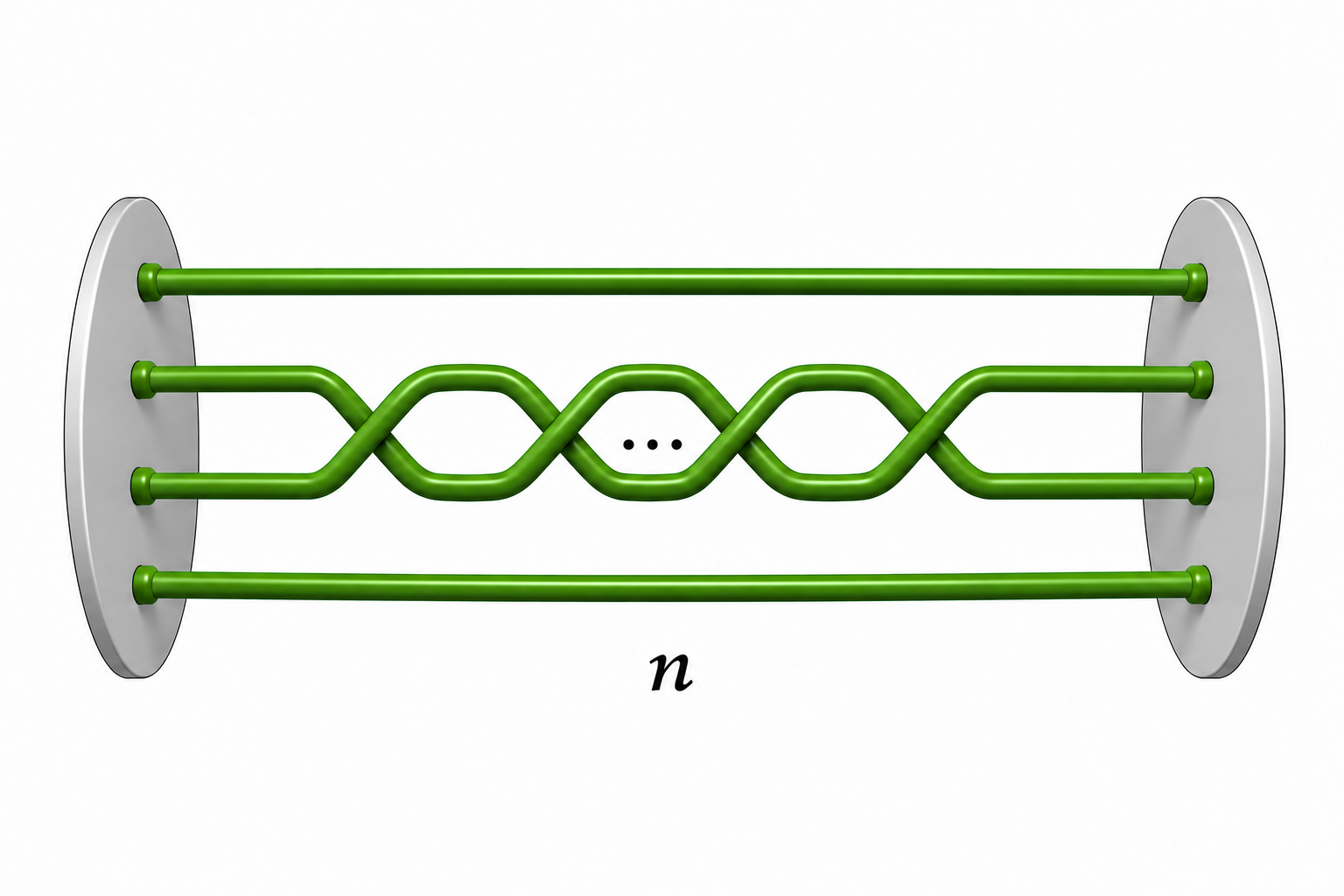}
}
&
The amplitude associated with two vertices($v_1$ \& $v_2$):
\[
P_{v_1v_2}\in\mathcal{H}^{\otimes 2}.
\]
\\
\hline

\end{tabular}

\caption{The graphical building blocks of an FRD and their TQFT. A
vertex combines states, a finger carries local braiding information, and a
propagator joins the states associated with two vertices.}
\label{tqft-tensor}
\end{table}

\vspace{-5mm} 
\paragraph{From an invariant to a knot name.}
The Explorer follows this same journey for every model. The user first
assembles a diagram by placing vertices, attaching fingers, drawing
propagators, and entering the signed twist parameters. Once the construction
is complete, the Explorer replaces each graphical element with its TQFT
counterpart and evaluates the resulting tensor network symbolically. Because
the currently supported FRDs are trees, the contraction can be organized from
the outer fingers toward the interior, allowing large diagrams to be built
from smaller calculations without introducing a closed cycle.

The final Laurent polynomial is then compared with a database of invariants
for classified knots. When the polynomial and its normalization agree with a
stored entry, the Explorer reports the corresponding knot as a database match.
This is an identification by the selected invariant, not a claim that the
polynomial uniquely determines the knot in every case: distinct knots can
occasionally share the same invariant. Where several entries agree, the result
should therefore be read as a set of candidates rather than a unique proof of
identity.

The Explorer also checks for {\it mirror images}. Mirroring reverses the order of the polynomial's coefficients, so the Explorer compares both the computed coefficient sequence and its reversal against the database, and reports when a match is obtained through the mirror relation.\footnote{For knots $9_{42}$ and 
$10_{71}$, the $r$-colored Jones polynomial with $r>2$ distinguishes from its mirror image, unlike the ordinary Jones, HOMFLY, and Kauffman/Akutsu--Wadati polynomials\cite{Ramadevi:1993hu}.} When the necessary data are available, the Explorer also displays the knot determinant as an additional invariant.

\subsection{Database: Knot Identification}

To identify knots, the Explorer compares the computed colored polynomial invariants with several curated polynomial databases. For knots with smaller crossing numbers (up to $10$ crossings), we use a database of colored invariants from the website \cite{Knotebook}, while the Adjoint polynomial ($r=2$) database from \cite{Singh2025} is used to identify knots with up to $13$ crossings. Whenever both the polynomial and its normalization coincide with a stored entry, the corresponding knot is reported as a database match. Once a match is found, the website displays the corresponding knot diagram using figures from \cite{KnotInfo}, providing a visual confirmation of the identification. These polynomial databases and knot diagrams will be expanded and updated regularly as new data become available.

\smallskip 

\textbf{How the results are checked.}
A successful symbolic calculation must also reproduce cases that are already
known. New constructions are therefore tested by \emph{reduction}: their
parameters are specialized until the diagram becomes a simpler, previously
verified model. For example, a three-vertex chain must reproduce the relevant
two-vertex result when one of its outer parts is reduced, and a branched
diagram must agree with a known chain when its additional branch is
specialized to the corresponding reduction value. These comparisons are made
term by term at the level of the Laurent polynomial. Only constructions that
pass the appropriate reduction tests are treated as verified.

\vspace{-5mm} 
\paragraph{Three ways to explore.}
The website presents the same mathematical pipeline through three levels of
construction, allowing a reader to begin with a familiar family and move
gradually toward a fully customizable FRD.

\vspace{-1mm} 
\begin{itemize}
\vspace{-2mm}
    \item[\bf (i)] 
\textbf{Pretzel diagrams.}
The pretzel interface is the most direct starting point. The current version
allows between three and seven fingers, each carrying a signed crossing count.
The user selects the representation color, enters the integers, and computes
the invariant. This mode is useful for reproducing standard examples and for
seeing immediately how a change in one twist region alters the polynomial.

\vspace{-2mm}
 \item[\bf (ii)] \textbf{Two-vertex trees.}
The next model introduces the essential idea of gluing local pieces. Two
vertices are connected by a propagator, and each vertex carries its own set of
fingers. Because the fingers may contain several marked braiding positions,
this model reaches diagrams that cannot be represented by the single-vertex
pretzel interface while remaining visually easy to follow.

\vspace{-2mm}
 \item[\bf (iii)] \textbf{Build your own.}
The free-form editor exposes the full tree-shaped construction. The user places
vertices on a canvas, attaches as many fingers as required, and draws
propagators between selected pairs of vertices. The direction in which a
propagator is drawn is purely visual; the engine records which vertices are
connected, not whether the line was drawn from left to right or from right to
left. A sequence of connected vertices forms a chain, while a vertex joined to
three or more neighboring vertices forms a branch. Both are supported, and
the interface displays the symbolic engine call generated by the diagram.

\end{itemize}

\vspace{-1mm} 
Closed cycles are deliberately excluded from the current implementation. A
cycle requires a different contraction strategy and does not share the same
tree factorization. Rather than applying an unsuitable formula, the Explorer
rejects such a diagram and makes the present boundary of the framework
explicit.

\begin{callout}{What the Explorer returns}
At present, the database provides SU$(2)$ Chern--Simons knot invariants as exact Laurent polynomials in $q$, computed for a chosen color (spin label) $r$. Each entry contains the following information:
\par\smallskip
\textbf{The invariant} --- the exact SU$(2)$ Chern--Simons knot invariant, expressed as a Laurent polynomial in \(q\), for the selected color $r$.
\par\smallskip
\textbf{The determinant} --- the sum of the absolute values of the coefficients of the invariant. For the ordinary Jones polynomial of an alternating knot, this quantity agrees with the classical knot determinant.
\par\smallskip
\textbf{The database match} ---the name of a knot with the same stored
invariant, together with a mirror indication when the match is obtained after reversing the powers of $q$.
\par\smallskip
\textcolor{blue}{\textbf{Future updates}} --- the database will be expanded to include SU$(N)$ Chern--Simons knot invariants for symmetric and selected mixed representations, the reformulated integer invariants inspired by the Chern--Simons/topological string correspondence, and many more computational and interactive database features.
\end{callout}
The Explorer is designed to let curiosity lead naturally to calculation. A reader can begin with a familiar pretzel knot, move to a two-vertex example,
and then build a new chain or branch from the same fundamental pieces. At each
stage, the diagram remains visible, the parameters remain explicit, and the
result can be reproduced.

To begin, choose a model, enter the signed twist values, select the
representation color, and run the calculation. The Explorer then follows the
same route that underlies the mathematics: from a visual FRD, to a TQFT tensor
network, to an exact invariant, and finally to a comparison with classified
knots. The result is more than a polynomial on a screen; it is a transparent
bridge between the shape of a knot and the quantum information it carries.

\medskip 
\begin{figure}[H]
    \centering
    \includegraphics[width=0.39\textwidth]{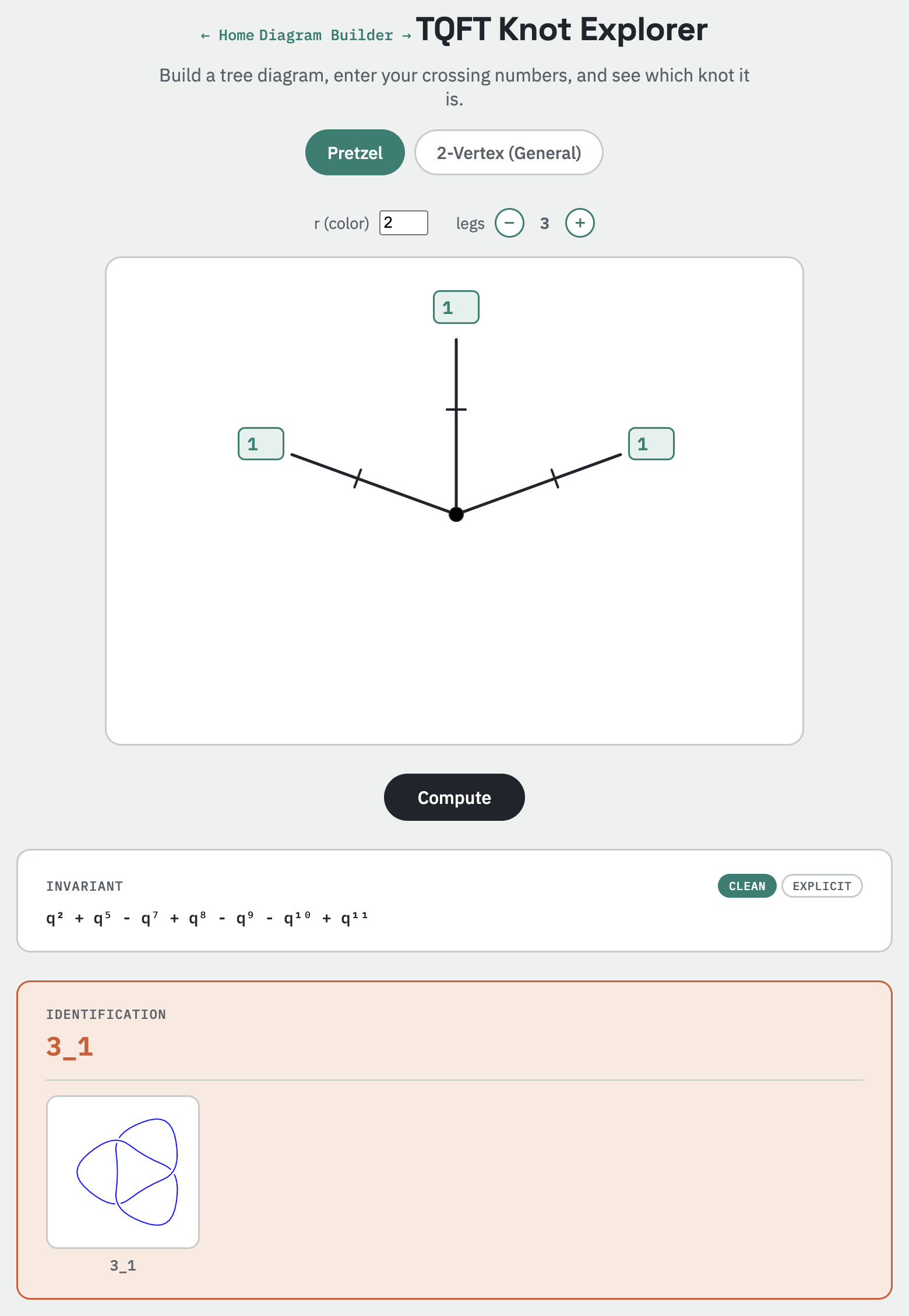}\qquad 
    \includegraphics[width=0.405\textwidth]{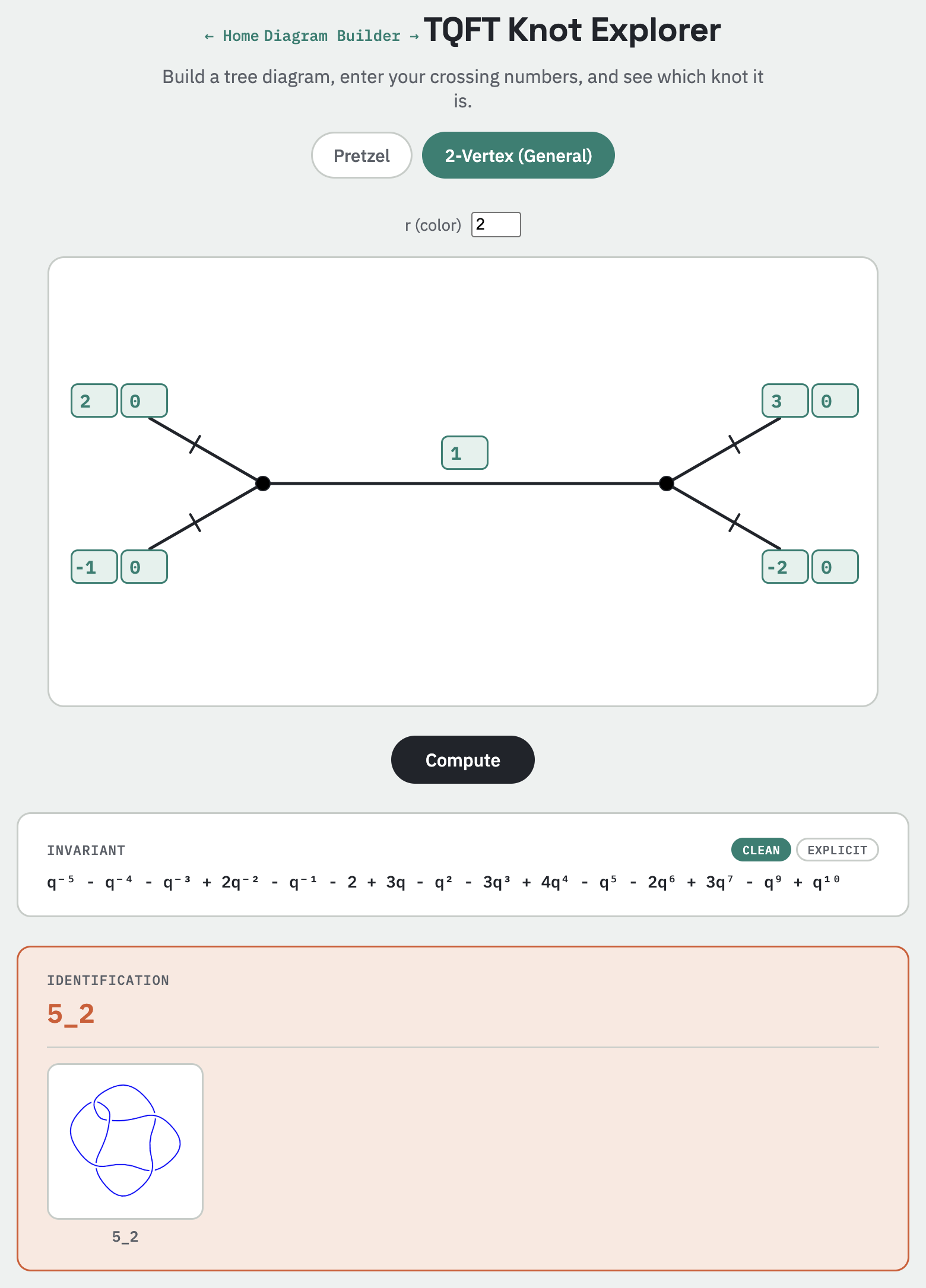}

    \vspace{2mm}

    \includegraphics[width=0.7\textwidth]{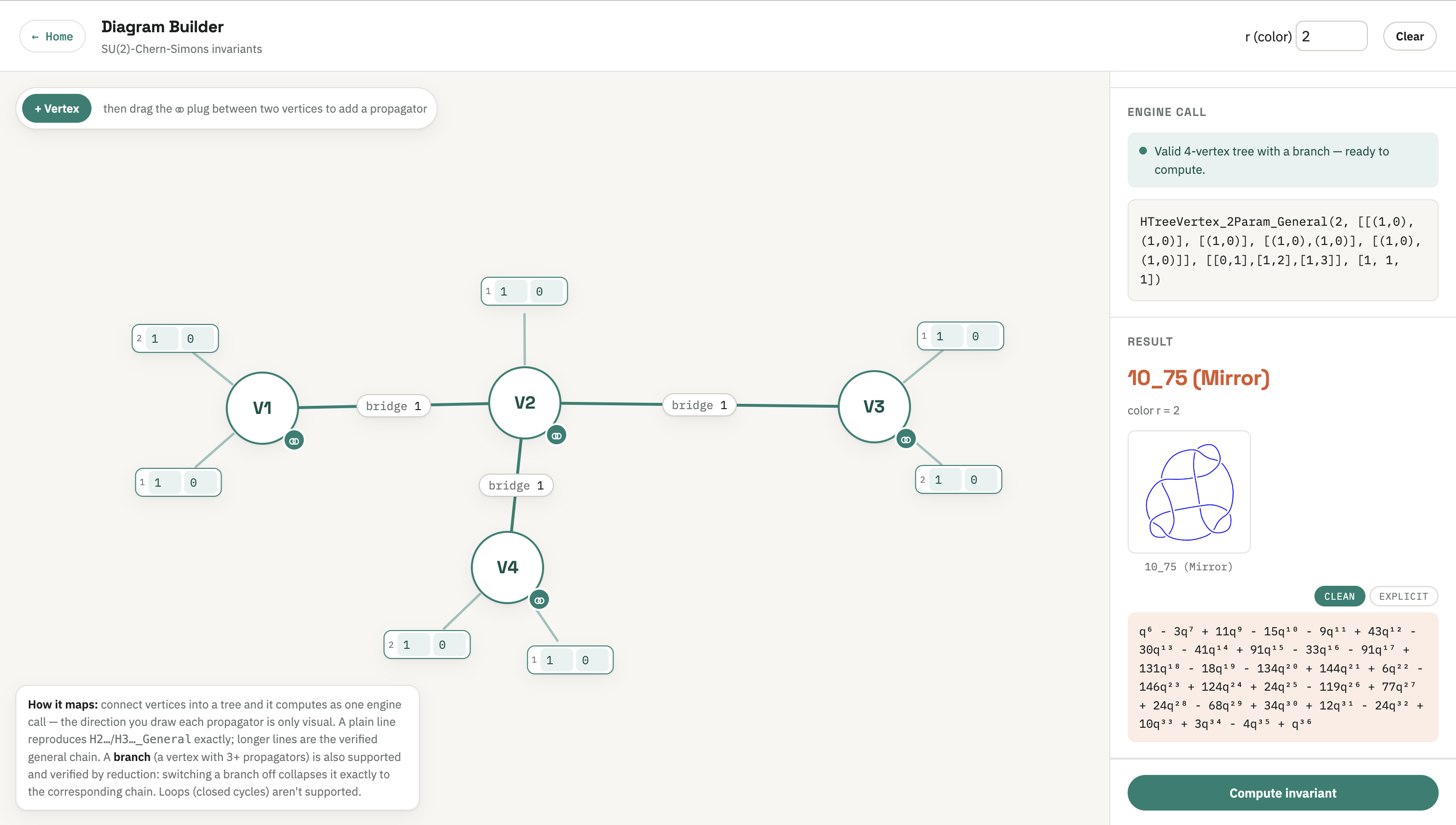}

    \caption{The web interface showing (left) the Pretzel explorer, (right) the 2-vertex explorer, and (bottom) the Diagram Builder.}
    \label{webinterface}
\end{figure}
\section{Efficient Algorithm and Search Procedure}
 In this section, we present our efficient algorithm and search procedure. We first introduce the notation and the conformal data of SU$(N)$-WZNW model used throughout our construction. Let $r$ denote a symmetric representation of SU$(N)$. We denote by $\mathrm{Qdim}(r)$ the quantum
dimension of the representation $r$, while $d(s,r)$ and $h(i,r)$ denote the quantum dimensions of the irreducible representations
$s\in r\otimes\bar{r}$ and $i\in r\otimes r$, respectively. The program involves two types of quantum $6j$-symbols, denoted by $a(i,j,r)$ and $e(i,j,r)$, together with two types of braiding eigenvalues: $\lambda^{+}_{i+1}$ for parallel strands and $\lambda^{-}_{s+1}$ for antiparallel strands.

The implementation currently available on the website computes colored Jones polynomials, corresponding to the SU$(2)$ case. In this case, the two channels are not distinguished: $a=e$, $d=h$, and, within our conventions, $\lambda^{+}=\lambda^{-}$. The extension to general SU$(N)$ requires a careful treatment of strand orientations and the associated fusion channels. This SU$(N)$ implementation will be added to the website in a forthcoming update\cite{RSS}. Further theoretical details will be presented in a companion work \cite{RSS2026}; see also Ref.~\cite{Nawata2013} for more context.
\subsection{Algorithmic Implementation}

The computation is implemented through a hybrid C++/Python pipeline. The
symbolic TQFT building blocks are evaluated in C++ using GiNaC, while the
enumeration of parameter families and knot identification are handled in
Python. The C++ routines are exposed to Python through \texttt{pybind11}.
The Python layer normalizes the resulting polynomial, converts it to a
coefficient representation, and compares it with a database of Jones and Adjoint ($r=2$) polynomials for knots up to $13$ crossings.
For an
internal state $i$, the one-parameter leg contribution is defined by
\[
P(i,m)
=
\sum_{s=0}^{r}
d(s,r)\,e(i,s,r)\,
\big(\lambda^-_{s+1}\big)^m
a(s,0,r).
\]
The two-parameter extension is
\vspace{-2mm} 
\[
P(i,m,n)
=
\sum_{s=0}^{r}
d(s,r)\,e(i,s,r)\,
\big(\lambda^-_{s+1}\big)^m
B_s(n),
\]
where
\[
B_s(n)=
\begin{cases}
a(s,0,r), & n=0,\\[4pt]
\displaystyle
\sum_{t=0}^{r}
d(t,r)\,a(s,t,r)
\big(\lambda^-_{t+1}\big)^n
a(t,0,r), & n\neq 0.
\end{cases}
\]
Thus, the two-parameter leg function reduces exactly to the
one-parameter leg function when $n=0$.

{\small
\begin{algorithm}[t]
\caption{\colorbox{lightgray}{Two-Vertex Propagator (independent indices)}}
\label{alg:two_vertex}
\begin{algorithmic}[1]
\Require Color $r$; leg data $\mathcal{A}$ for the first vertex;
bridge exponent $b$; leg data $\mathcal{B}$ for the second vertex.
\Ensure The independent-index two-vertex invariant $H_{2V}^{\mathrm{gen}}$.
 
\State $v_A \gets |\mathcal{A}|+1$, \quad $v_B \gets |\mathcal{B}|+1$
\State $H \gets 0$
 
\For{$t_1 = 0,\ldots,r$}
    \State
    \vspace{-2mm}
    \[
    V_A(t_1)
    \gets
    \frac{
    h(t_1,r)\prod_{(m,n)\in\mathcal{A}} P(t_1,m,n)
    }{
    e(t_1,0,r)^{\,v_A-2}
    }
    \]
    \State $R \gets 0$ \Comment{row accumulator for this $t_1$}
    \For{$t_2 = 0,\ldots,r$}
        \State
        \[
        V_B(t_2)
        \gets
        \frac{
        h(t_2,r)\prod_{(m,n)\in\mathcal{B}} P(t_2,m,n)
        }{
        e(t_2,0,r)^{\,v_B-2}
        }
        \]
        \State
        \[
        R \gets R +
        V_A(t_1)\,\operatorname{Prop}(t_1,t_2,b)\,V_B(t_2)
        \]
    \EndFor
    \State $H \gets H + \operatorname{Simplify}(R)$
    \Comment{simplify each row before accumulating}
\EndFor
 
\State \Return $\operatorname{Simplify}\!\left((-1)^r\operatorname{Qdim}(r)\,H\right)$
\end{algorithmic}
\end{algorithm}
}

The generalized propagator $\operatorname{Prop}(t_1,t_2,b)$ is defined in the
following subsection. Algorithm~\ref{alg:two_vertex} specializes in two
independent ways. Forcing the two vertices to share a single internal state,
$t_1=t_2=t$, collapses the double sum to a single sum and replaces
$\operatorname{Prop}$ by the bare bridge factor
$(\lambda^-_{t+1})^{b}/h(t,r)$, recovering the shared-index invariant
$H_{2V}$. Independently, restricting every leg to the form $(m,0)$ gives
$P(t,m,0)=P(t,m)$, so the one-parameter model is the $n=0$ restriction of the
two-parameter one. The two specializations may be applied separately or
together.

\subsection{Independent-Index Generalization}

To allow the two vertices to carry independent internal states, the bridge
factor is replaced by the generalized propagator
\[
\operatorname{Prop}(t_1,t_2,b)
=
\sum_{s=0}^{r}
d(s,r)\,a(t_1,s,r)
\big(\lambda^-_{s+1}\big)^b
a(s,t_2,r).
\]
The independent-index two-vertex invariant is therefore
\[
H_{2V}^{\mathrm{gen}}
=
(-1)^r\operatorname{Qdim}(r)
\sum_{t_1=0}^{r}\sum_{t_2=0}^{r}
V_A(t_1)\,
\operatorname{Prop}(t_1,t_2,b)\,
V_B(t_2).
\]
In the implementation, the inner sum over $t_2$ is simplified row by row
before being added to the total sum. This avoids a single large symbolic
simplification at the end and substantially improves performance.

{\small 
\begin{algorithm}[H]
\caption{\colorbox{lightgray}{Enumeration and Knot Identification}}
\label{alg:scan}
\begin{algorithmic}[1]
\Require Maximum crossing budget $S_{\max}$; bridge exponent $b$;
knot polynomial database $\mathcal{D}$.
\Ensure List of parameter configurations whose polynomials match known knots.

\State Initialize an empty set $\mathcal{S}$ of evaluated canonical forms

\For{$S = 4, \ldots, S_{\max}$}
    \State Build pair pool $\mathcal{P} \gets \bigl\{(m,n) : m \neq 0,\;
    |m|+|n| = c,\; c = 1,\ldots,S-3 \bigr\}$,
    enumerating both signs of $m$ and $n$
    \For{each 4-multiset $\{(m_1,n_1),(m_2,n_2),(m_3,n_3),(m_4,n_4)\}$
    drawn from $\mathcal{P}$ via combinations with replacement
    such that $\sum_{i}(|m_i|+|n_i|) = S$}
        \For{each partition of the $4$-multiset into ordered vertices
        $(\mathcal{A},\mathcal{B})$ with $|\mathcal{A}|=|\mathcal{B}|=2$}
            \State $\mathcal{A} \gets \operatorname{sort}(\mathcal{A})$,\quad
                   $\mathcal{B} \gets \operatorname{sort}(\mathcal{B})$
            \Comment{canonicalize legs within each vertex}
            \State $C \gets \operatorname{sort}\bigl(\mathcal{A},\mathcal{B}\bigr)$
            \Comment{order the two vertices}
            \If{$C \in \mathcal{S}$}
                \State \textbf{continue}
            \EndIf
            \State Add $C$ to $\mathcal{S}$
            \State Compute $H_{2V}^{\mathrm{gen}}(r,\,\mathcal{A},\,b,\,\mathcal{B})$
            using Algorithm~\ref{alg:two_vertex}
            \State Normalize, convert to a coefficient tuple, and look it up
            in $\mathcal{D}$ (as above)
            \If{a match is found}
                \State Record the knot name and $(\mathcal{A},b,\mathcal{B})$
            \EndIf
        \EndFor
    \EndFor
\EndFor
\end{algorithmic}
\end{algorithm}
}

\vspace{-3mm} 
The canonical form in Algorithm~\ref{alg:scan} implements the symmetry
of the two-vertex construction, which acts at the level of vertices. Two operations leave the invariant unchanged: permuting the legs \emph{within} a single vertex, since each vertex factor is a
product over its legs; and exchanging the two vertices as complete blocks,
$(\mathcal{A},\mathcal{B})\mapsto(\mathcal{B},\mathcal{A})$, since the
generalized propagator is symmetric in its indices,
$\operatorname{Prop}(t_1,t_2,b)=\operatorname{Prop}(t_2,t_1,b)$. Moving a leg
\emph{across} the bridge from one vertex to the other is \emph{not} a symmetry:
the two legs of a vertex are multiplied together and attached to a single
propagator index, so redistributing them changes the invariant. A given
multiset of four pairs therefore admits several inequivalent two-vertex
diagrams, one per distinct $2{+}2$ partition up to the vertex swap, and each is
scanned separately. Exchanging the two entries \emph{within} a pair,
$(m_i,n_i)\mapsto(n_i,m_i)$, is likewise not a symmetry, since the two
parameters enter distinct summations. The bridge exponent is treated as a fixed
scan parameter.

Finally, knot identification is performed by hashing the coefficient tuple
of each normalized polynomial. A second hash table stores reversed
coefficient tuples, allowing mirror images to be detected in constant time.

\subsection{Computational Scalability of Colored Jones Invariants}
\begin{figure}[H]
    \centering
    \begin{adjustbox}{center}
        \includegraphics[width=18cm]{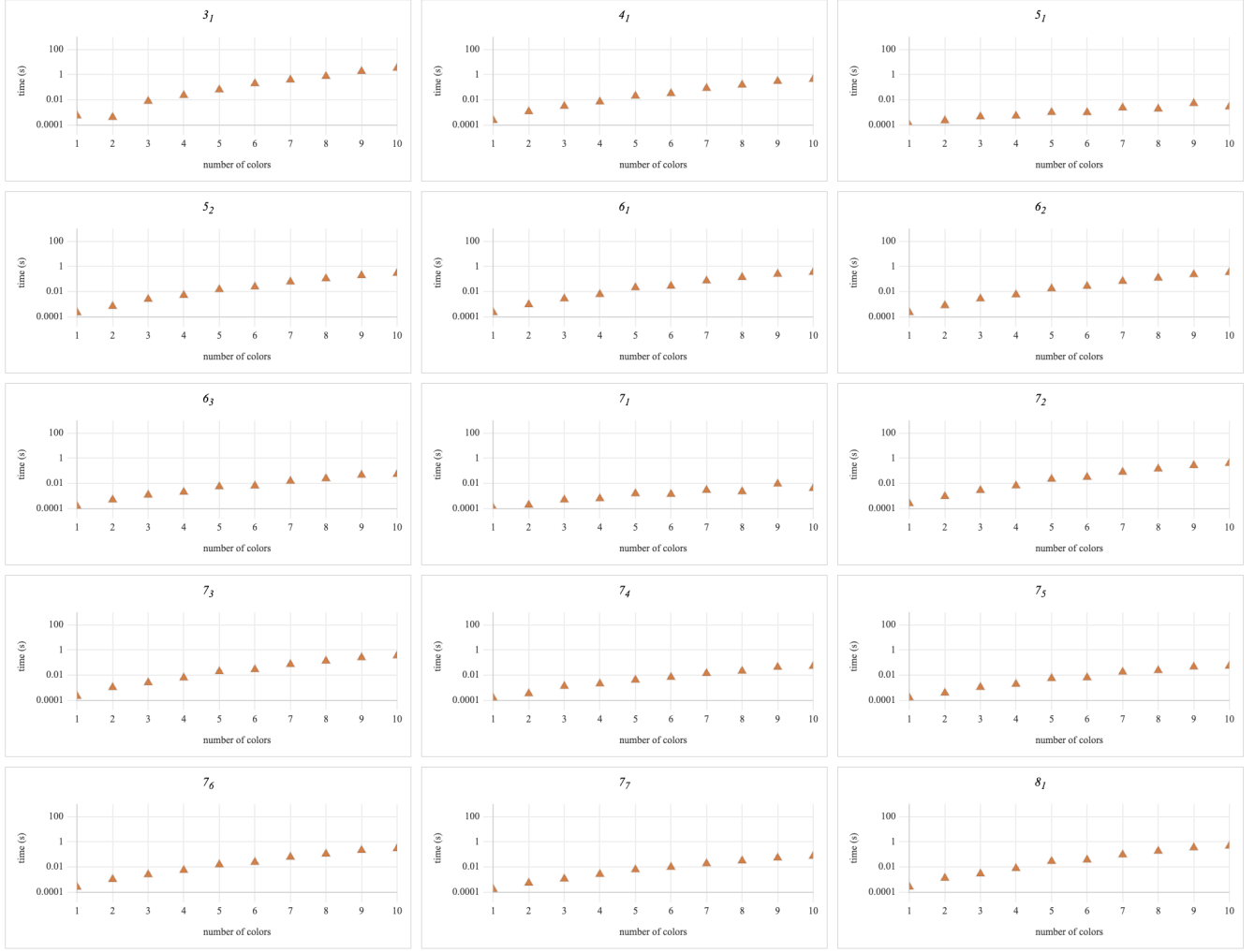}
    \end{adjustbox}
    \caption{Running time of the single-vertex pretzel algorithm
    (\texttt{HPretzelK}) implemented in C++ (GiNaC) for the first 15
    arborescent knots up to 8 crossings, at representation parameter
    $r = 1, \ldots, 10$.}
    \label{2vertex}
\end{figure}

\begin{figure}[H]
    \centering
    \begin{adjustbox}{center}
        \includegraphics[width=18cm]{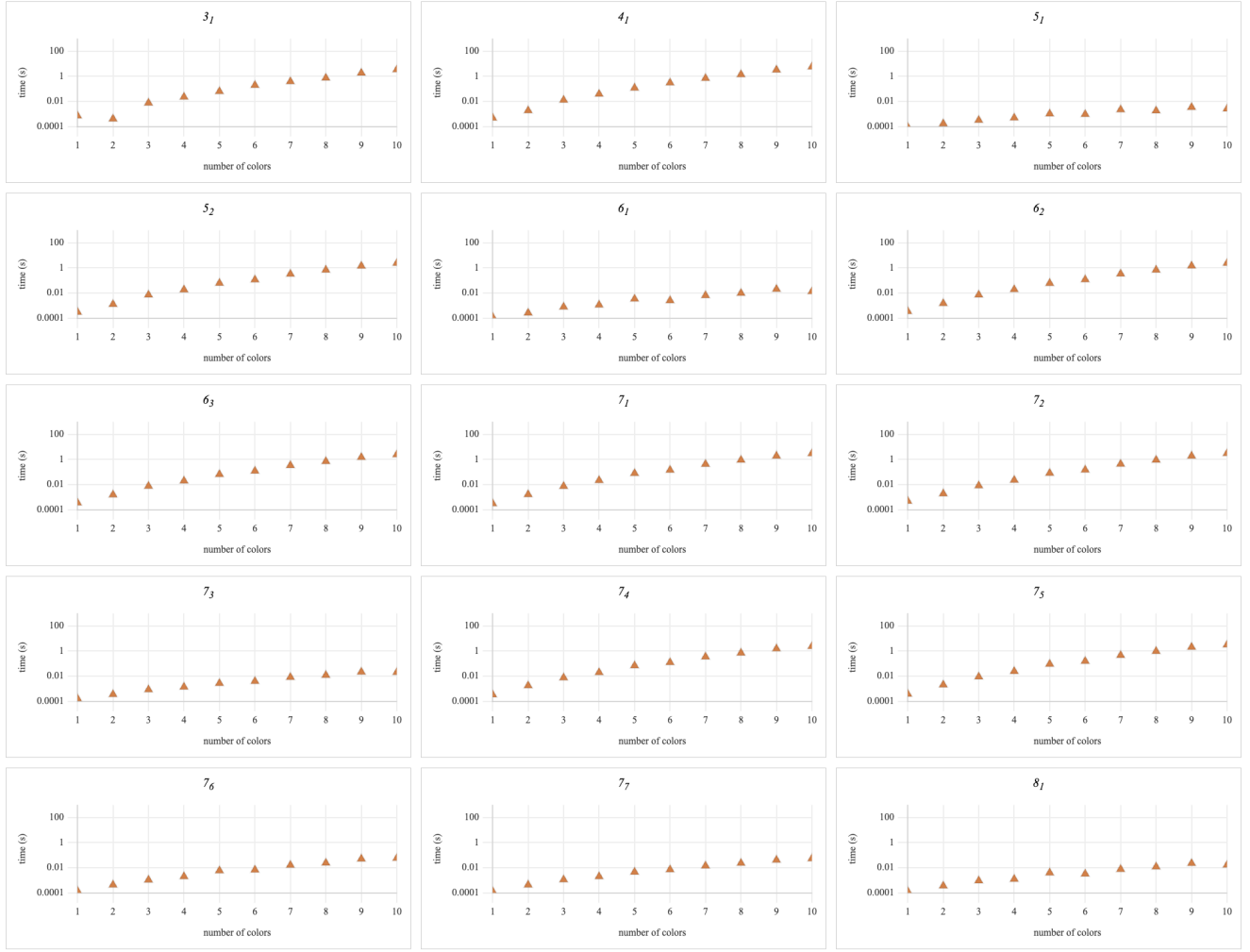}
    \end{adjustbox}
    \caption{Running time of the 2-vertex 2-parameter tree propagator
    algorithm (\texttt{H2VertexPropagatorK\_2Param}) implemented in C++
    (GiNaC) for the first 15 arborescent knots up to 8 crossings, at
    representation parameter $r = 1, \ldots, 10$.}
    \label{testpretzel}
\end{figure}

\begin{figure}[H]
    \centering
    \begin{adjustbox}{center}
        \includegraphics[width=18cm]{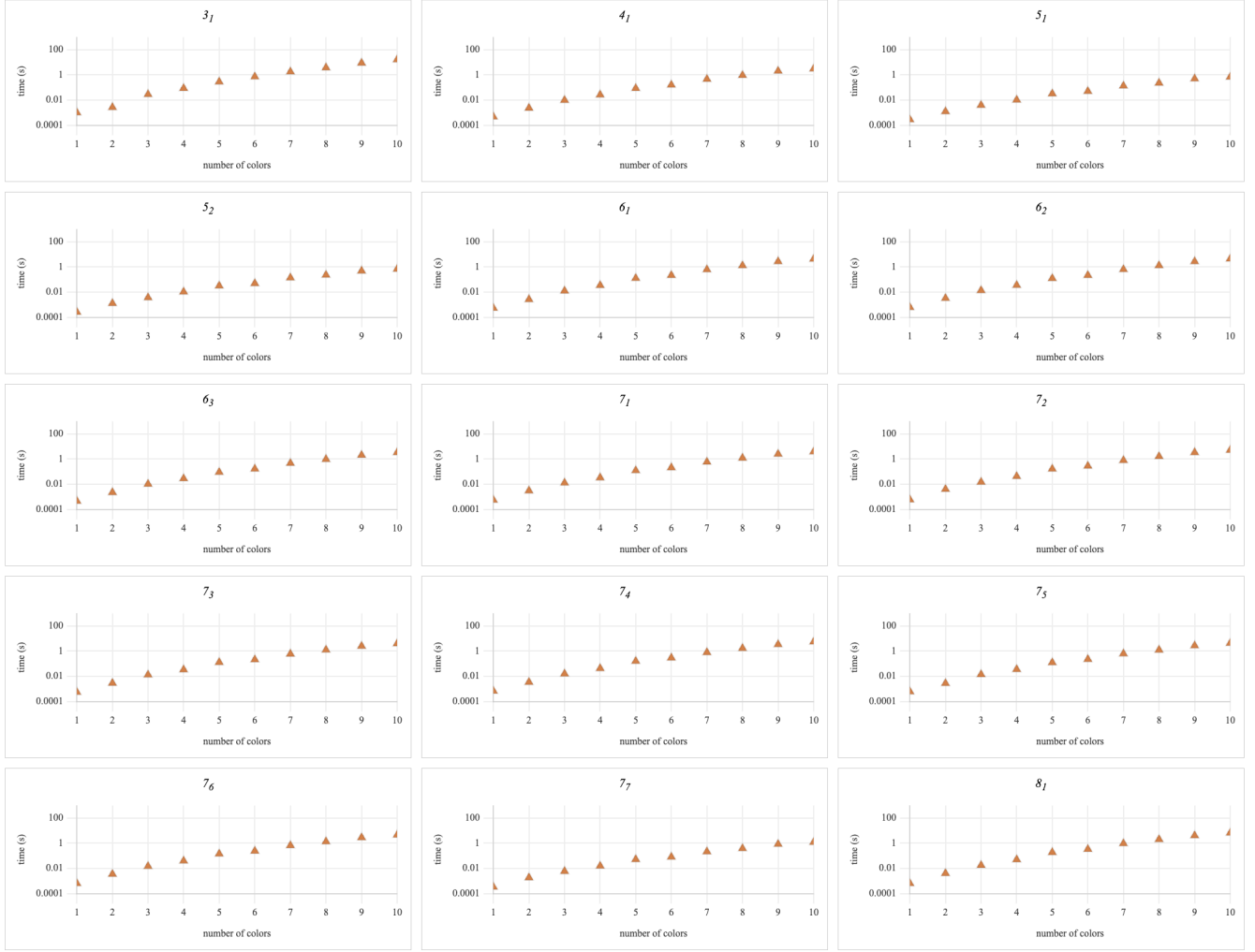}
    \end{adjustbox}
    \caption{Running time of the 2-vertex 2-parameter general tree propagator
    algorithm \newline (\texttt{H2VertexPropagatorK\_General\_2Param}) implemented in C++
    (GiNaC) for the first 15 arborescent knots up to 8 crossings, at
    representation parameter $r = 1, \ldots, 10$.}{\label{22vertex}}
\end{figure}

\begin{figure}[H]
    \centering
    \begin{adjustbox}{center}
        \includegraphics[width=18cm]{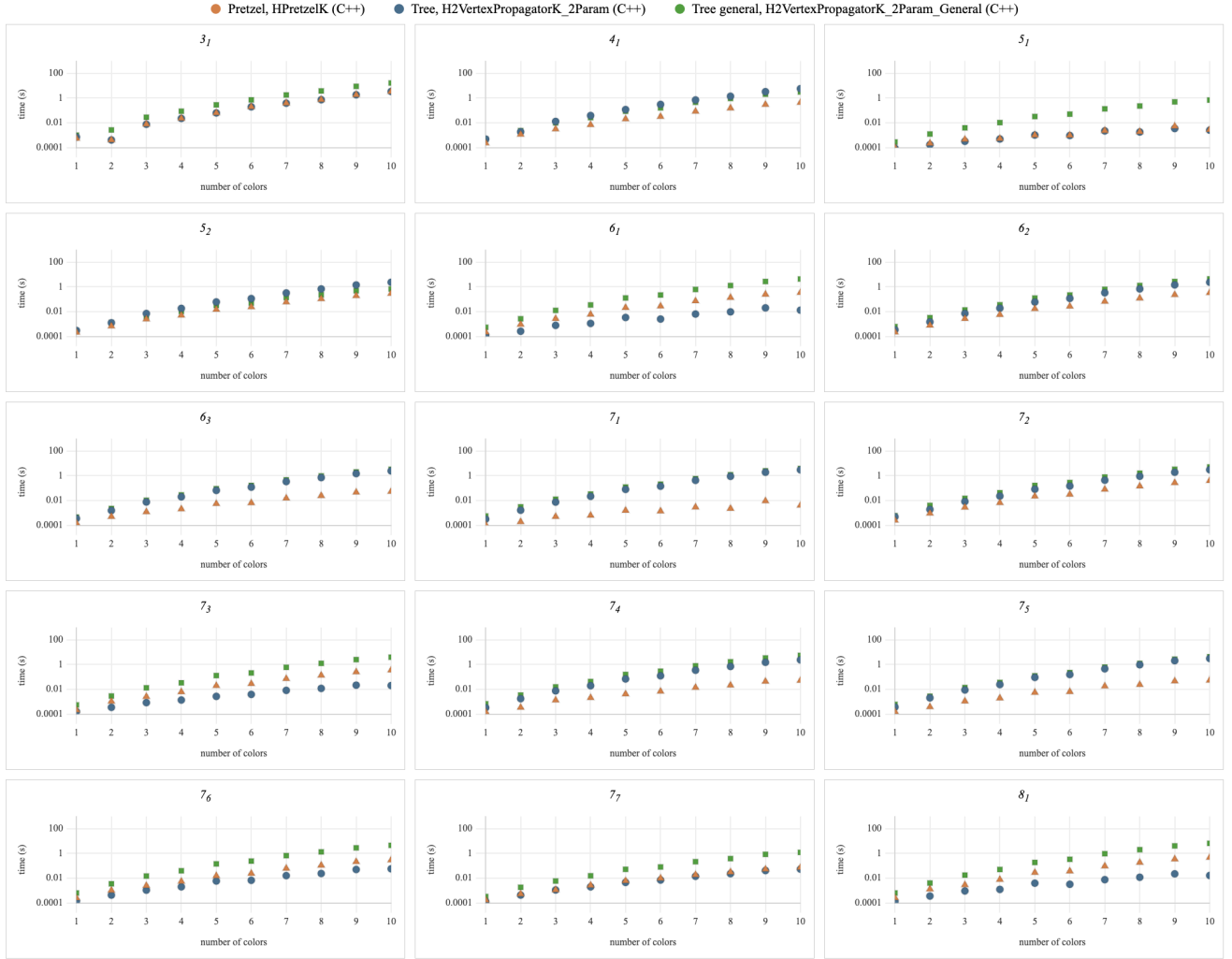}
    \end{adjustbox}
\caption{Comparison of the single-vertex pretzel algorithm and the
    2-vertex 2-parameter tree propagator algorithm, both implemented in
    C++ (GiNaC), at representation parameter $r = 1, \ldots, 10$.}
    \label{pretzel-vs-tree}
\end{figure}

We analyze the computational scalability of our algorithm for colored Jones invariants on the proposed evolution families, namely pretzel knots, one-vertex evolution families, and two-vertex evolution families; see ~\cite{RSS2026} for their definitions and constructions. The benchmark results, presented in Figures~\ref{2vertex}, \ref{testpretzel}, \ref{22vertex} and \ref{pretzel-vs-tree}, show that the proposed method maintains a stable and efficient running time even as the structural complexity of the knots increases.

\smallskip 
For a systematic evaluation, we compare our implementation with the CJP algorithm from the \texttt{KnotTheory} package~\cite{KnotAtlas} and with the braid-walk-based approach of Ref.~\cite{HajijLevitt2018}. The optimized Python implementation developed in this work reduces unnecessary computational overhead by taking advantage of the recursive/evolution structure of the knot families. As a result, in most tested cases, our method achieves faster computation than existing approaches.

\smallskip 
The significance of these benchmarks is twofold. First, they confirm that the proposed evolution-based algorithm is scalable for increasingly complex knot families. Second, they show that the method provides a practical computational framework for exploring colored Jones invariants beyond the range that is typically accessible using more general-purpose algorithms. This makes the approach particularly well-suited for large-scale experimental investigations in quantum knot invariants.

\bigskip 
\bigskip 

\section{Discussion and Conclusion}\label{dis-con}

We have developed a systematic and user-friendly web-based platform\cite{RSS} for studying higher-rank invariants of knots represented by Feynman ribbon diagrams (FRDs), also known as arborescent diagrams (see Fig.\ref{fig:knot-10-93}) in the literature\cite{Con,Caudron,BS}. This broad class contains several important families of knots, including pretzel knots and two-bridge knots. The Explorer allows a user to draw an arbitrary FRD tree, identify the corresponding knot, investigate its structural classification, and compute its colored knot polynomials within a unified computational framework.

A central feature of the Explorer is its automated knot-identification pipeline. After a user submits a tree diagram to the FRD Explorer, the system first computes a collection of strong knot invariants.\footnote{Here, \emph{strong invariants} refers to a collection of invariants capable of distinguishing most knots in the database. Nevertheless, a few exceptional cases remain. For example, colored Jones polynomials alone cannot, in general, distinguish mutant knots. These remaining ambiguities can be eliminated in future refinements by extending the database to include higher-rank colored invariants.}
These invariants are then compared with the entries in the associated knot database to identify the corresponding knot whenever a match is available. Once the identification is completed, the user can efficiently compute higher-rank invariants  without performing the full symbolic calculation manually.

The Explorer represents a significant step toward a unified computational framework for quantum knot invariants. The present restriction to $13$-crossing knots reflects only the available knot database, while the computational engine itself remains fully general for admissible FRD inputs. Future developments~\cite{RSS} will extend the platform to higher-crossing knots and SU$(N)$ gauge groups, establishing it as a versatile resource for modern quantum topology.

Another central contribution of this work is the two-vertex FRD structure, which we call the \emph{evolution diagram}; see Fig.~\ref{fig:FRD}. This compact structure reveals how families of FRD-like knots are related to one another. It therefore provides a systematic classification of all such knots with up to $10$ crossings. More importantly, the same construction extends naturally to knots with higher crossing numbers. To our knowledge, an FRD-based classification beyond $10$ crossings has not previously been presented in this unified form~\cite{MMMRSS2017}.

The value of the framework is not limited to classification. Our benchmarks against the \texttt{KnotTheory} colored Jones algorithm and the braid-walk method~\cite{KnotAtlas,HajijLevitt2018} show that it is faster in most of the tested cases. Thus, the proposed method offers an efficient and scalable approach to both knot classification and invariant computation. The accompanying Explorer makes the framework accessible for further investigation. It allows users to construct, modify, and analyze their own FRD diagrams. In this way, existing classifications can be refined and extended toward a more complete classification of arborescent knots.

Beyond knot classification, the invariants introduced here provide a foundation for future investigations in several emerging research directions, including the knot--quiver correspondence~\cite{Kucharski2019,StosicWedrich2021,Chauhan:2025iwt}, machine-learning methods for knot invariants~\cite{CravenJejjalaKar2021,CravenHughesJejjalaKar2024,HughesJejjalaRamadeviRoySingh2025}, and quantum computing~\cite{Kauffman2002,Passante2009,Rasetti:2006hzb}.  Overall, this work establishes the FRD framework as a link between structural classification, efficient computation, and broader applications of knot invariants.

\section{Appendix}{\label{Appendix1}}
Using the evolution family of two-vertex Feynman ribbon diagrams, $\Gamma(A_1;p_{12};A_2)$ (see Fig. \ref{fig:FRD2} and the discussion around it), we obtain a complete classification of the FRD types for all knots with up to $10$ crossings.
\begingroup
\scriptsize
\setlength{\tabcolsep}{2pt}
\renewcommand{\arraystretch}{1.05}
\begin{longtable}{|p{0.02\textwidth}|p{0.33\textwidth}|p{0.055\textwidth}|p{0.33\textwidth}|p{0.055\textwidth}|}
\hline
\rowcolor{lightgray}
\centering\textbf{ }
& \centering\textbf{FRD Parameters}
& \centering\textbf{Knot}
& \centering\textbf{FRD Parameters}
& \centering\arraybackslash\textbf{Knot} \\ \hline
\endfirsthead

\hline
\centering\textbf{}
& \centering\textbf{FRD Parameters}
& \centering\textbf{Knot}
& \centering\textbf{FRD Parameters}
& \centering\arraybackslash\textbf{Knot} \\ \hline
\endhead

\hline
\endfoot

\multirow{7}{*}{
}
& $([(-1, -1), (-1, 0)], 1, [(1, 0), (1, 0)])$ & $3_1$
& $([(-3, 1), (-1, 0)], 1, [(1, 0), (1, 0)])$ & $7_1$ \\ \cline{2-5}

& $([(-1, 0), (-1, 0)], 1, [(-1, 1), (1, 0)])$ & $4_1$
& $([(-2, 1), (-1, 0)], 1, [(1, -1), (1, 0)])$ & $7_2$ \\ \cline{2-5}

& $([(-1, 0), (-1, 1)], 1, [(1, 0), (1, 0)])$ & $5_1$
& $([(-1, 0), (-1, 1)], 1, [(1, 0), (2, -1)])$ & $7_3$ \\ \cline{2-5}

& $([(-1, 0), (-1, 0)], 1, [(1, -1), (1, 0)])$ & $5_2$
& $([(-1, 1), (1, 0)], 1, [(1, 0), (2, -1)])$ & $7_4$ \\ \cline{2-5}

& $([(-1, 0), (-1, 1)], 1, [(1, -1), (1, 0)])$ & $6_1$
& $([(-1, 0), (-1, 0)], 1, [(1, -1), (2, -1)])$ & $7_5$ \\ \cline{2-5}

& $([(-1, 1), (1, 0)], 1, [(1, -1), (1, 0)])$ & $6_2$
& $([(-2, 1), (1, 0)], 1, [(1, -1), (1, 0)])$ & $7_6$ \\ \cline{2-5}

& $([(1, -1), (1, 0)], 1, [(1, 0), (1, 0)])$ & $6_3$
& $([(1, -1), (1, 0)], 1, [(1, -1), (1, 0)])$ & $7_7$ \\ \hline

\end{longtable}
\endgroup

\begingroup
\scriptsize
\setlength{\tabcolsep}{2pt}
\renewcommand{\arraystretch}{1.05}
\begin{longtable}{|p{0.02\textwidth}|p{0.33\textwidth}|p{0.055\textwidth}|p{0.33\textwidth}|p{0.055\textwidth}|}
 \hline
 \rowcolor{lightgray}
\centering\textbf{ }
& \centering\textbf{FRD Parameters}
& \centering\textbf{Knot}
& \centering\textbf{FRD Parameters}
& \centering\arraybackslash\textbf{Knot} \\ \hline
\endfirsthead

\hline
\rowcolor{lightgray}
\centering\textbf{}
& \centering\textbf{FRD Parameters}
& \centering\textbf{Knot}
& \centering\textbf{FRD Parameters}
& \centering\arraybackslash\textbf{Knot} \\ \hline
\endhead

\hline
\endfoot

\multirow{32}{*}{
}
& $([(-3, 1), (-1, 0)], 1, [(1, -1), (1, 0)])$ & $8_1$
& $([(-1, 1), (1, 0)], 1, [(1, 0), (2, -3)])$ & $9_{13}$ \\ \cline{2-5}

& $([(-1, 3), (1, 0)], 1, [(1, -1), (1, 0)])$ & $8_2$
& $([(1, -1), (1, 0)], 1, [(1, 0), (3, -1)])$ & $9_{14}$ \\ \cline{2-5}

& $([(-1, 0), (-1, 1)], 1, [(1, 0), (3, -1)])$ & $8_3$
& $([(-1, 0), (-1, 0)], 1, [(1, -1), (3, -2)])$ & $9_{15}$ \\ \cline{2-5}

& $([(-1, 1), (1, 0)], 1, [(1, 0), (3, -1)])$ & $8_4$
& $([(1, 2), (1, 0)], -1, [(-2, 1), (-2, 1)])$ & $9_{16}$ \\ \cline{2-5}

& $([(1, -1), (-1, 0)], -1, [(-1, 1), (-2, 1)])$ & $8_5$
& $([(1, 0), (1, -1)], -1, [(-1, 2), (-1, 2)])$ & $9_{17}$ \\ \cline{2-5}

& $([(-1, 0), (-1, 1)], 1, [(1, -1), (2, -1)])$ & $8_6$
& $([(-1, 0), (-1, 0)], 1, [(1, -1), (2, -3)])$ & $9_{18}$ \\ \cline{2-5}

& $([(-1, 0), (-1, 0)], -1, [(-1, 0), (-3, 1)])$ & $8_7$
& $([(-1, 0), (-1, 0)], -1, [(-1, 0), (-3, 2)])$ & $9_{19}$ \\ \cline{2-5}

& $([(1, 0), (1, -1)], -1, [(-1, 2), (-1, 1)])$ & $8_8$
& $([(-2, 2), (1, 0)], 1, [(1, 0), (2, -1)])$ & $9_{20}$ \\ \cline{2-5}

& $([(-1, 1), (1, 0)], 1, [(1, -3), (1, 0)])$ & $8_9$
& $([(1, -1), (1, 2)], 1, [(1, 0), (2, -1)])$ & $9_{21}$ \\ \cline{2-5}

& $([(-1, 1), (-2, 1)], -1, [(-1, 0), (-1, 0)])$ & $8_{10}$
& $([(-1, 1), (-1, 0)], -1, [(-1, 1), (-2, 1)])$ & $9_{22}$ \\ \cline{2-5}

& $([(-3, 1), (1, 0)], 1, [(1, -1), (1, 0)])$ & $8_{11}$
& $([(2, -1), (-1, 0)], -1, [(-1, 0), (-2, 2)])$ & $9_{23}$ \\ \cline{2-5}

& $([(-1, 0), (-1, 0)], 1, [(1, -1), (2, -2)])$ & $8_{12}$
& $([(-1, 1), (-1, -4)], -1, [(-1, 0), (-1, 0)])$ & $9_{24}$ \\ \cline{2-5}

& $([(1, -1), (1, 0)], 1, [(1, 0), (2, -1)])$ & $8_{13}$
& $([(-1, 1), (-2, 2)], -1, [(-1, 0), (-1, 0)])$ & $9_{25}$ \\ \cline{2-5}

& $([(-1, 0), (-1, 0)], -1, [(-1, 0), (-2, 2)])$ & $8_{14}$
& $([(1, -3), (1, 0)], 1, [(1, -1), (1, 0)])$ & $9_{26}$ \\ \cline{2-5}

& $([(-1, 2), (-1, 1)], -1, [(-1, 0), (-1, 0)])$ & $8_{15}$
& $([(1, -2), (1, 2)], 1, [(1, -1), (1, 0)])$ & $9_{27}$ \\ \cline{2-5}

& $([(-1, 1), (-1, 2)], 1, [(-1, 1), (-1, 2)])$ & $8_{16}$
& $([(-1, 2), (-1, -3)], -1, [(-1, 0), (-1, 0)])$ & $9_{28}$ \\ \cline{2-5}

& $([(-2, 0), (-1, 2)], 1, [(-2, 1), (-1, 1)])$ & $8_{17}$
& $([(-1, 1), (-1, 2)], 1, [(-1, 2), (-1, 2)])$ & $9_{29}$ \\ \cline{2-5}

& $([(-2, 0), (-2, 1)], 1, [(1, 0), (1, 0)])$ & $8_{19}$
& $([(-1, 2), (-1, 1)], -1, [(-1, 1), (-1, 0)])$ & $9_{30}$ \\ \cline{2-5}

& $([(-1, 1), (-1, 2)], 1, [(1, 0), (1, 0)])$ & $8_{20}$
& $([(-1, 2), (-1, 0)], -1, [(-1, 2), (-1, 0)])$ & $9_{31}$ \\ \cline{2-5}

& $([(1, -1), (-1, 2)], -1, [(-1, 0), (-1, 0)])$ & $8_{21}$
& $([(-2, 0), (-2, 1)], 1, [(-2, 0), (-2, 2)])$ & $9_{32}$ \\ \cline{2-5}

& $([(5, -1), (1, 0)], -1, [(-1, 0), (-1, 0)])$ & $9_1$
& $([(2, 0), (1, -2)], -1, [(2, -2), (1, -1)])$ & $9_{33}$ \\ \cline{2-5}

& $([(-4, 1), (-1, 0)], 1, [(1, -1), (1, 0)])$ & $9_2$
& $([(1, -1), (-1, 0)], -1, [(-2, 1), (-2, 1)])$ & $9_{35}$ \\ \cline{2-5}

& $([(-3, 1), (-1, 0)], 1, [(1, 0), (2, -1)])$ & $9_3$
& $([(1, -1), (-1, 0)], -1, [(-1, 1), (-2, 2)])$ & $9_{36}$ \\ \cline{2-5}

& $([(-2, 1), (-1, 0)], 1, [(1, 0), (3, -1)])$ & $9_4$
& $([(-1, 2), (-2, 1)], -1, [(-1, 0), (-1, 0)])$ & $9_{37}$ \\ \cline{2-5}

& $([(-1, 1), (1, 0)], 1, [(1, 0), (4, -1)])$ & $9_5$
& $([(-2, 0), (-1, 2)], 1, [(-2, 1), (-2, 1)])$ & $9_{38}$ \\ \cline{2-5}

& $([(-1, 0), (-1, 0)], 1, [(1, -1), (4, -1)])$ & $9_6$
& $([(-2, 0), (-2, 1)], 1, [(1, -1), (1, 0)])$ & $9_{42}$ \\ \cline{2-5}

& $([(-2, 1), (-1, 0)], 1, [(1, -1), (2, -1)])$ & $9_7$
& $([(1, -1), (-2, 1)], -1, [(-1, 1), (-1, 0)])$ & $9_{43}$ \\ \cline{2-5}

& $([(2, -1), (1, 0)], -1, [(-1, 2), (-1, 1)])$ & $9_8$
& $([(-1, 1), (-1, 2)], 1, [(1, -1), (1, 0)])$ & $9_{44}$ \\ \cline{2-5}

& $([(-1, 0), (-1, 0)], 1, [(2, -1), (3, -1)])$ & $9_9$
& $([(1, -1), (-1, 2)], -1, [(-1, 1), (-1, 0)])$ & $9_{45}$ \\ \cline{2-5}

& $([(-1, 0), (-1, 1)], 1, [(2, -1), (2, -1)])$ & $9_{10}$
& $([(-2, 1), (-1, 2)], 1, [(1, 0), (1, 0)])$ & $9_{46}$ \\ \cline{2-5}

& $([(-1, 0), (-1, 0)], 1, [(1, -4), (1, -1)])$ & $9_{11}$
& $([(2, -1), (-1, 2)], -1, [(-1, 0), (-1, 0)])$ & $9_{48}$ \\ \cline{2-5}

& $([(-4, 1), (1, 0)], 1, [(1, -1), (1, 0)])$ & $9_{12}$
& & \\ \hline

\end{longtable}
\endgroup

\begingroup
\scriptsize
\setlength{\tabcolsep}{2pt}
\renewcommand{\arraystretch}{1.05}
\begin{longtable}{|p{0.02\textwidth}|p{0.33\textwidth}|p{0.055\textwidth}|p{0.33\textwidth}|p{0.055\textwidth}|} \hline
\rowcolor{lightgray}
\centering\textbf{ }
& \centering\textbf{FRD Parameters}
& \centering\textbf{Knot}
& \centering\textbf{FRD Parameters}
& \centering\arraybackslash\textbf{Knot} \\ \hline
\endfirsthead

\hline
\rowcolor{lightgray}
\centering\textbf{}
& \centering\textbf{FRD Parameters}
& \centering\textbf{Knot}
& \centering\textbf{FRD Parameters}
& \centering\arraybackslash\textbf{Knot} \\ \hline
\endhead

\hline
\endfoot

\multirow{65}{*}{
}
& $([(-5, 1), (-1, 0)], 1, [(1, -1), (1, 0)])$ & $10_1$
& $([(-1, 3), (-1, 2)], -1, [(-1, 0), (-1, 0)])$ & $10_{66}$ \\ \cline{2-5}

& $([(-1, 5), (1, 0)], 1, [(1, -1), (1, 0)])$ & $10_2$
& $([(-1, 0), (-1, 0)], -1, [(-2, 2), (-2, 1)])$ & $10_{67}$ \\ \cline{2-5}

& $([(-3, 1), (-1, 0)], 1, [(1, 0), (3, -1)])$ & $10_3$
& $([(-1, 1), (-1, 0)], -1, [(-2, 1), (-2, 1)])$ & $10_{68}$ \\ \cline{2-5}

& $([(-1, 1), (1, 0)], 1, [(1, 0), (5, -1)])$ & $10_4$
& $([(-1, 2), (-1, 2)], -1, [(-1, 1), (-1, 0)])$ & $10_{69}$ \\ \cline{2-5}

& $([(-1, 0), (-1, 0)], -1, [(-1, 0), (-5, 1)])$ & $10_5$
& $([(1, 2), (1, 0)], -1, [(-2, 2), (-2, 1)])$ & $10_{70}$ \\ \cline{2-5}

& $([(-1, 0), (-1, 1)], 1, [(1, -1), (4, -1)])$ & $10_6$
& $([(1, 0), (1, 0)], 1, [(1, 3), (2, -2)])$ & $10_{71}$ \\ \cline{2-5}

& $([(-5, 1), (1, 0)], 1, [(1, -1), (1, 0)])$ & $10_7$
& $([(-1, 1), (-1, 0)], -1, [(-1, 1), (-1, -4)])$ & $10_{72}$ \\ \cline{2-5}

& $([(-1, 2), (1, 0)], 1, [(1, 0), (4, -1)])$ & $10_8$
& $([(-1, 2), (-1, -3)], -1, [(-1, 1), (-1, 0)])$ & $10_{73}$ \\ \cline{2-5}

& $([(-1, 1), (1, 0)], 1, [(1, -5), (1, 0)])$ & $10_9$
& $([(-1, 0), (-1, 0)], -1, [(-1, -4), (-2, 1)])$ & $10_{74}$ \\ \cline{2-5}

& $([(1, -1), (1, 0)], 1, [(1, 0), (4, -1)])$ & $10_{10}$
& $([(1, 2), (1, -1)], -1, [(-1, 2), (-1, 2)])$ & $10_{75}$ \\ \cline{2-5}

& $([(-1, 0), (-1, 1)], 1, [(2, -1), (3, -1)])$ & $10_{11}$
& $([(1, 0), (1, -2)], -1, [(-2, 1), (-2, 1)])$ & $10_{76}$ \\ \cline{2-5}

& $([(1, 0), (1, -1)], -1, [(-1, 2), (-3, 1)])$ & $10_{12}$
& $([(1, 0), (1, -2)], -1, [(-1, 2), (-2, 1)])$ & $10_{77}$ \\ \cline{2-5}

& $([(-1, 0), (-1, 0)], 1, [(1, -1), (2, -4)])$ & $10_{13}$
& $([(1, 0), (1, -2)], -1, [(-1, 2), (-1, 2)])$ & $10_{78}$ \\ \cline{2-5}

& $([(-1, 0), (-1, 0)], -1, [(-1, 0), (-2, 4)])$ & $10_{14}$
& $([(-1, -4), (-1, 1)], 1, [(-1, 1), (2, -1)])$ & $10_{79}$ \\ \cline{2-5}

& $([(-1, 0), (-1, 1)], 1, [(1, -4), (1, -1)])$ & $10_{15}$
& $([(-1, -3), (-1, 2)], 1, [(-1, 1), (2, -1)])$ & $10_{80}$ \\ \cline{2-5}

& $([(-1, 0), (-1, 0)], 1, [(1, -4), (2, -1)])$ & $10_{16}$
& $([(-1, -3), (-1, 2)], 1, [(-1, 1), (1, -2)])$ & $10_{81}$ \\ \cline{2-5}

& $([(-1, 2), (1, 0)], 1, [(1, -4), (1, 0)])$ & $10_{17}$
& $([(-2, 0), (-1, 4)], 1, [(-2, 1), (-1, 1)])$ & $10_{82}$ \\ \cline{2-5}

& $([(1, -1), (1, 2)], 1, [(1, 0), (3, -1)])$ & $10_{18}$
& $([(-2, 0), (-2, 1)], 1, [(-2, 0), (-2, 3)])$ & $10_{83}$ \\ \cline{2-5}

& $([(1, 0), (2, -1)], 1, [(1, 0), (3, -1)])$ & $10_{19}$
& $([(3, -2), (1, -1)], -1, [(1, -1), (1, -2)])$ & $10_{84}$ \\ \cline{2-5}

& $([(-3, 1), (-1, 0)], 1, [(1, -1), (2, -1)])$ & $10_{20}$
& $([(-1, 1), (-1, 2)], 1, [(-1, 1), (-1, 4)])$ & $10_{85}$ \\ \cline{2-5}

& $([(2, -1), (1, 0)], -1, [(-1, 2), (-2, 1)])$ & $10_{21}$
& $([(-2, 0), (-1, 2)], 1, [(-2, 3), (-1, 1)])$ & $10_{86}$ \\ \cline{2-5}

& $([(-1, 1), (1, 0)], 1, [(1, 0), (3, -3)])$ & $10_{22}$
& $([(3, 0), (1, -1)], -1, [(3, -2), (1, -1)])$ & $10_{87}$ \\ \cline{2-5}

& $([(-1, 0), (-1, 0)], -1, [(-1, 0), (-3, 3)])$ & $10_{23}$
& $([(-2, -3), (-1, 1)], 1, [(-2, 0), (-2, 2)])$ & $10_{88}$ \\ \cline{2-5}

& $([(-1, 0), (-1, 1)], 1, [(1, -1), (2, -3)])$ & $10_{24}$
& $([(2, 0), (2, -2)], -1, [(2, 0), (2, -2)])$ & $10_{89}$ \\ \cline{2-5}

& $([(-3, 3), (1, 0)], 1, [(1, -1), (1, 0)])$ & $10_{25}$
& $([(-2, 0), (-2, 1)], 1, [(-1, 2), (-1, 3)])$ & $10_{90}$ \\ \cline{2-5}

& $([(1, 0), (1, 4)], 1, [(1, 0), (2, -1)])$ & $10_{26}$
& $([(-2, 0), (-2, 1)], 1, [(-2, 1), (-1, 3)])$ & $10_{91}$ \\ \cline{2-5}

& $([(-1, 1), (-1, 0)], -1, [(-1, 0), (-2, 3)])$ & $10_{27}$
& $([(2, 0), (1, -2)], -1, [(2, -2), (1, -2)])$ & $10_{92}$ \\ \cline{2-5}

& $([(1, 0), (1, -1)], -1, [(-1, 3), (-1, 2)])$ & $10_{28}$
& $([(-1, 1), (-1, 2)], 1, [(-1, 2), (-1, 3)])$ & $10_{93}$ \\ \cline{2-5}

& $([(-3, 2), (1, 0)], 1, [(1, 0), (2, -1)])$ & $10_{29}$
& $([(-3, 0), (-1, 1)], 1, [(-3, 1), (-1, 2)])$ & $10_{94}$ \\ \cline{2-5}

& $([(1, -3), (1, 2)], 1, [(1, -1), (1, 0)])$ & $10_{30}$
& $([(-2, 0), (-2, 1)], 1, [(-2, 2), (-1, 2)])$ & $10_{95}$ \\ \cline{2-5}

& $([(-1, 0), (-1, -2)], -1, [(-1, 0), (-3, 2)])$ & $10_{31}$
& $([(-2, 0), (-2, 2)], 1, [(-1, 2), (-1, 2)])$ & $10_{96}$ \\ \cline{2-5}

& $([(-1, 0), (-2, 2)], -1, [(-1, 0), (-2, 1)])$ & $10_{32}$
& $([(2, 0), (2, -2)], -1, [(2, -1), (2, -1)])$ & $10_{97}$ \\ \cline{2-5}

& $([(1, -3), (1, 0)], 1, [(1, 0), (2, -1)])$ & $10_{33}$
& $([(-2, 1), (-1, 2)], 1, [(-1, 2), (-1, 2)])$ & $10_{98}$ \\ \cline{2-5}

& $([(3, -1), (1, 0)], -1, [(-1, 2), (-1, 1)])$ & $10_{34}$
& $([(-2, 1), (-1, 2)], 1, [(-2, 1), (-1, 2)])$ & $10_{99}$ \\ \cline{2-5}

& $([(-1, 0), (-1, 0)], 1, [(1, -1), (4, -2)])$ & $10_{35}$
& $([(4, -1), (1, -1)], -1, [(-1, 0), (-1, 0)])$ & $10_{124}$ \\ \cline{2-5}

& $([(-1, 0), (-1, 0)], -1, [(-1, 0), (-4, 2)])$ & $10_{36}$
& $([(1, -1), (-4, 1)], -1, [(-1, 0), (-1, 0)])$ & $10_{125}$ \\ \cline{2-5}

& $([(-1, 0), (-1, 1)], 1, [(1, -1), (3, -2)])$ & $10_{37}$
& $([(-1, 1), (-1, 4)], 1, [(1, 0), (1, 0)])$ & $10_{126}$ \\ \cline{2-5}

& $([(2, -1), (-1, 0)], -1, [(-1, 0), (-3, 2)])$ & $10_{38}$
& $([(1, -1), (-1, 4)], -1, [(-1, 0), (-1, 0)])$ & $10_{127}$ \\ \cline{2-5}

& $([(1, 0), (1, -1)], -1, [(-1, 2), (-2, 2)])$ & $10_{39}$
& $([(-2, 0), (-2, 1)], 1, [(1, 0), (2, -1)])$ & $10_{128}$ \\ \cline{2-5}

& $([(1, -1), (1, 0)], 1, [(1, 2), (2, -2)])$ & $10_{40}$
& $([(-1, 1), (-1, 2)], 1, [(1, 0), (2, -1)])$ & $10_{129}$ \\ \cline{2-5}

& $([(-2, 2), (1, 0)], 1, [(1, 0), (2, -2)])$ & $10_{41}$
& $([(1, -1), (-2, 1)], -1, [(-1, 0), (-2, 1)])$ & $10_{130}$ \\ \cline{2-5}

& $([(-1, 2), (-1, 0)], -1, [(-1, 0), (-2, 2)])$ & $10_{42}$
& $([(1, -1), (-1, 2)], -1, [(-1, 0), (-2, 1)])$ & $10_{131}$ \\ \cline{2-5}

& $([(-1, 0), (-1, 0)], 1, [(1, -2), (3, 3)])$ & $10_{43}$
& $([(-3, 2), (-1, 1)], 1, [(1, 0), (1, 0)])$ & $10_{132}$ \\ \cline{2-5}

& $([(1, -1), (1, 0)], 1, [(1, 0), (3, 3)])$ & $10_{44}$
& $([(1, -1), (-3, 2)], -1, [(-1, 0), (-1, 0)])$ & $10_{133}$ \\ \cline{2-5}

& $([(1, -2), (1, 0)], 1, [(1, 0), (2, 3)])$ & $10_{45}$
& $([(1, -1), (-2, 1)], -1, [(-1, 1), (-1, -2)])$ & $10_{134}$ \\ \cline{2-5}

& $([(1, -1), (-1, 0)], -1, [(-1, 1), (-4, 1)])$ & $10_{46}$
& $([(3, -2), (-1, 1)], -1, [(-1, 0), (-1, 0)])$ & $10_{135}$ \\ \cline{2-5}

& $([(-1, 1), (-4, 1)], -1, [(-1, 0), (-1, 0)])$ & $10_{47}$
& $([(-2, 0), (-2, 2)], 1, [(1, -1), (1, 0)])$ & $10_{136}$ \\ \cline{2-5}

& $([(1, -1), (-1, 0)], -1, [(-1, 4), (-1, 1)])$ & $10_{48}$
& $([(1, -1), (-2, 2)], -1, [(-1, 1), (-1, 0)])$ & $10_{137}$ \\ \cline{2-5}

& $([(-1, 4), (-1, 1)], -1, [(-1, 0), (-1, 0)])$ & $10_{49}$
& $([(2, -2), (-1, 1)], -1, [(-1, 1), (-1, 0)])$ & $10_{138}$ \\ \cline{2-5}

& $([(1, -1), (-1, 0)], -1, [(-1, 1), (-2, 3)])$ & $10_{50}$
& $([(-3, 0), (-3, 1)], 1, [(1, 0), (1, 0)])$ & $10_{139}$ \\ \cline{2-5}

& $([(-1, 1), (-2, 3)], -1, [(-1, 0), (-1, 0)])$ & $10_{51}$
& $([(-3, 1), (-1, 2)], 1, [(1, 0), (1, 0)])$ & $10_{140}$ \\ \cline{2-5}

& $([(-1, 1), (-2, 1)], -1, [(-1, 0), (-2, 1)])$ & $10_{52}$
& $([(1, -3), (-1, 2)], -1, [(-1, 0), (-1, 0)])$ & $10_{141}$ \\ \cline{2-5}

& $([(-1, 2), (-1, 1)], -1, [(-1, 0), (-2, 1)])$ & $10_{53}$
& $([(-2, 1), (-1, 3)], 1, [(1, 0), (1, 0)])$ & $10_{142}$ \\ \cline{2-5}

& $([(1, -1), (-1, 0)], -1, [(-1, 1), (-3, 2)])$ & $10_{54}$
& $([(-1, 2), (-1, 3)], 1, [(1, 0), (1, 0)])$ & $10_{143}$ \\ \cline{2-5}

& $([(-1, 1), (-3, 2)], -1, [(-1, 0), (-1, 0)])$ & $10_{55}$
& $([(2, -1), (-1, 3)], -1, [(-1, 0), (-1, 0)])$ & $10_{144}$ \\ \cline{2-5}

& $([(-1, 1), (-1, -2)], -1, [(-1, 1), (-2, 1)])$ & $10_{56}$
& $([(-2, 1), (-2, 1)], 1, [(1, -1), (1, 0)])$ & $10_{145}$ \\ \cline{2-5}

& $([(-1, 2), (-1, 1)], -1, [(-1, 1), (-1, -2)])$ & $10_{57}$
& $([(-1, 2), (-1, 2)], 1, [(1, -1), (1, 0)])$ & $10_{146}$ \\ \cline{2-5}

& $([(2, -1), (-1, 0)], -1, [(-1, 1), (-2, 2)])$ & $10_{58}$
& $([(1, -2), (-1, 2)], -1, [(-1, 1), (-1, 0)])$ & $10_{147}$ \\ \cline{2-5}

& $([(-1, 1), (-1, 0)], -1, [(-1, 1), (-2, 2)])$ & $10_{59}$
& $([(-2, 0), (-2, 1)], 1, [(-1, 1), (2, -1)])$ & $10_{148}$ \\ \cline{2-5}

& $([(-1, 1), (-1, 0)], -1, [(-1, 1), (-2, -3)])$ & $10_{60}$
& $([(-1, 1), (-1, 2)], 1, [(-1, 1), (2, -1)])$ & $10_{149}$ \\ \cline{2-5}

& $([(1, -1), (-1, 0)], -1, [(-2, 1), (-3, 1)])$ & $10_{61}$
& $([(-2, 0), (-2, 1)], 1, [(-1, 1), (1, -2)])$ & $10_{150}$ \\ \cline{2-5}

& $([(-1, 0), (-1, 0)], -1, [(-2, 1), (-3, 1)])$ & $10_{62}$
& $([(-1, 1), (-1, 2)], 1, [(-1, 1), (1, -2)])$ & $10_{151}$ \\ \cline{2-5}

& $([(-1, 2), (-3, 1)], -1, [(-1, 0), (-1, 0)])$ & $10_{63}$
& $([(-1, 1), (2, -1)], 1, [(-1, 1), (2, -1)])$ & $10_{152}$ \\ \cline{2-5}

& $([(1, -1), (-1, 0)], -1, [(-1, 3), (-2, 1)])$ & $10_{64}$
& $([(-1, 1), (1, -2)], 1, [(-1, 1), (2, -1)])$ & $10_{153}$ \\ \cline{2-5}

& $([(-1, 3), (-2, 1)], -1, [(-1, 0), (-1, 0)])$ & $10_{65}$
& $([(-1, 1), (1, -2)], 1, [(-1, 1), (1, -2)])$ & $10_{154}$ \\ \hline

\end{longtable}
\endgroup

\bigskip 

\end{document}